\documentclass[10pt]{article}

\usepackage{graphicx}
\usepackage{graphicx,psfrag}
\usepackage{amsmath,amssymb,amsthm}
 
\def \vs {\vskip 0.2cm}

\def \n {\noindent}

\def \ds {\displaystyle }

\def \la {{\langle}}

\def \ra {{\rangle}}

\def \CA {{\cal A}}

\def \CV {{\cal V}}
\def \CE {{ \cal E}}

\def \CW {{\cal W}}

\def \CL {{\cal L}}

\def \CC {{\cal C}}

\def \CF {{\cal F}}

\def \CG {{\cal G}}

\def \O {{\Omega}}

\def \E {{\bf E}}

\def \bZ {{\mathbb Z}}

\def \bE {{\mathbb E}}

\def \bR {{\mathbb R}}

\def \bP {{\mathbb P}}

\def \bN {{\mathbb N}}

\def \g {\gamma }

\def \z {\zeta}

\def \D {\Delta}

\def \Tr {{\hbox {\,Tr\,}}}

\def \Th {\Theta}

\def \bE {{\mathbb E}}

\def \BY {{\bf Y}}

\def \a {\alpha}
 
\def \b {\beta} 

\def \s { \sigma}

\def \l {\lambda}

\def \L {\Lambda}

\def \d {\delta}

\def \vep {\varepsilon}

\def \ep {\epsilon}

\def \k {\kappa}
\def \u {\upsilon}

\def \U {{\Upsilon}}

\def \G {\Gamma}

\def \CA {{\cal A}}

\def \CD {{\cal D}}

\def \CW {{\cal W}}

\def \CV {{\cal V}} 

\def \CC {{\cal C}}

\def \CT {{\cal T}}

\def \CH {{\cal H}}

\def \CE {{\cal E}}

\def \CC {{\cal C}}

\def \CN {{\cal N}}

\def \CX {{\cal X}}

\def \CY {{\cal Y}}

\def \CK{{\cal K}}

\def \CM {{\cal M}}

\def \CP {{\cal P}}

\def \CL {{\cal L}}

\def \CS {{\cal S}}

\def \fD {{\mathfrak D}}

\def \fT {{\mathfrak T}}

\def \fN {{ \mathfrak N}}

\def \fL {{\mathfrak L}}

\def \fP {{\mathfrak P}}

\def \fe {{\mathfrak e}}

\def \RM {{\rm M}}

\def \rmA{{\rm A}}
\def \rmB{{\rm B}}
\def \rmC{{\rm C}}

\def \rma{{\rm a}}
\def \rmb{{\rm b}}
\def \rmc{{\rm c}}

\def \rV {{\rm V}}
\def \rE {{\rm E}}

\def \rW{{\rm W}}
\def \rL {{\rm L}}

\def \rT {{\rm T}}

\def \ri {{\rm i}}
\def \rii {{\rm ii}}
\def \riii {{\rm iii}}

\def \rA {{\rm A}}

\def \rG {{\rm G}}

\def \rc {{\rm color}}

\def \rg {{\rm g}}

\def \rcon {{\rm conn}}

\def \rcum{{ \rm Cum}}

\def \rD {{\rm D}}

\def \rV {{\rm V}}

\def \rv {{\rm v}}

\begin{document}
\title{Limit theorems for 
walks and 
triangles on  Erd\H os-R\'enyi random graphs with large interaction radius
}

\author{O. Khorunzhiy\\ Universit\'e de Versailles - Saint-Quentin \\45, Avenue des Etats-Unis, 78035 Versailles, FRANCE\\
{\it e-mail:} oleksiy.khorunzhiy@uvsq.fr}
\maketitle
\begin{abstract}
We study 
 cumulants of numbers of $q$-step 
 walks on  
Erd\H os-R\'enyi   random graphs with distance-dependent edge probability 
in the limit when the number of vertices $N$, concentration $c$,   and interaction radius $R$ tend  to infinity.
These cumulants can be associated with   a formal   cumulant expansion of the 
free energy  of matrix models of exponential random graphs
widely known in mathematical and theoretical physics.

We show that 
in three different asymptotic regimes, 
the limiting values of \mbox{$k$-th} cumulants  $\CF_k^{(q)}$ exist and 
can be  associated with one or another family of tree-type diagrams, in dependence on the asymptotic behavior of parameters
$cR/N$ for $q$-step non-closed walks and $c^2R/N^2$ for 3-step closed walks, respectively. In some cases, 
we obtain explicit expressions for $\CF_k^{(q)}$ with the help of a modified Pr\"ufer  codification.

These results 
allow us to 
prove Limit Theorems for the number of non-closed walks  and for the number of triangles in corresponding ensembles of large random graphs. 
We indicate an asymptotic threshold that separates 
the  normal probability distribution and  the Poisson one for the number of triangles in random graphs. 
We
show that 
in  the random graph  ensemble 
that  we consider
the average vertex degree 
can be bounded from above while the total  
number of triangles infinitely increases,
thus rigorously solving a graph collapse problem known in applications.

  \end{abstract}



\section{Introduction}

Random matrix theory and random graphs theory are closely related,
in particular  
 by means of random adjacency matrices. 
Let us consider a family  $
\CA_n= \{ a_{ij}^{(n)}\}_{ 1\le i < j \le n}  
$
of  Bernoulli random variables determined on the same probability space
such that
$$
a_{ij}^{(n)} = 
\begin{cases}
 1,  &  \text{with probability $p_n(i,j)$}  , \\
0,  \quad   & 
\text {with probability $1-p_n(i,j)$
}, 
\end{cases}
\quad 1\le i  < j\le n. 
\eqno (1.1)
$$
A real symmetric random matrix $A_{n}$ with matrix elements $a^{(n)}_{ij}$
over the main diagonal and zeros  on the diagonal can be considered  as the adjacency matrix of  a simple non-oriented loop-less graph on  $n$ vertices. 
In this context, the values  $ p_n(i,j)$ can be  regarded as  edge probabilities of a random graph ensemble determined by (1.1), 
i.e.  
the probability for the edge $e(i,j)$ to be present in the random graph.

In the case when  $\CA_n$ is a family of jointly independent random variables such that the edge probabilities $p_n= p_n(i,j)$  do not depend on $i$ and $j$, 
one gets an
 ensemble of random graphs 
$\{A_{n}^{(p_n)}\}$ that in the limit of infinite $n$  
is asymptotically close to the Erd\H os-R\'enyi ensemble of random graphs
\cite{B}.
This ensemble $\{A_{n}^{(p_n)}\}$ is often referred  by itself as to the Erd\H os-R\'enyi (or Erd\H os-R\'enyi-Gilbert) ensemble of random graphs $\rG(n,p_n)$. 
Asymptotic properties of  the Erd\H os-R\'enyi random graphs in the limit of large dimension $n\to\infty$ are thoroughly investigated in vaste  number of papers.
 A large variety of properties and phenomena has been observed 
  in dependence of the limiting behavior
of the edge probability $p_n$ as $n$ tending to infinity. 
One of the results states  that  
the number of triangles in the Erd\H os-R\'enyi random graphs with $p_n= c/n$
converges in the limit  $n\to\infty$
to a random variable $\nu$
that follows the Poisson probability distribution
 \cite{B}, 
$$
T_n= \#\{\hbox{triangles in $\rG(n,c/n)$}\} 
\stackrel{\CL}{\to}\nu, \ \nu \sim  \CP(c^3/6), \  \ n\to\infty.
\eqno (1.2)
$$ 
Convergence (1.2) have been further studied in a  number of papers (see, for example, \cite{CRF,GHN} and references therein).

The use of adjacency matrices in subgraph counting can be fairly effective. In particular,
one can associate  variable $T_n$ (1.2) with the trace $X_n^{(3)} = \Tr A^3_n$ 
that can be also regarded as  the  total number of closed three-step 
walks over 
 the  graph determined \mbox{by $A_n=A_n^{(p_n)}$.}
One can define also the total number of two-step non-closed walks by 
 $Y^{(2)}_n= \sum_{i,j} (A^2_n)_{ij}$.

 In paper \cite{K-08}, asymptotic properties of the cumulant expansion of 
 variables 
  $$
Z_n^{(X^{(q)})}(g)= \log 
\bE_{{ER}} \exp \left\{ - g 
 X^{(q)}_n \right\},
\qquad X^{(q)}_n = \Tr A^q_n = \sum_{i} \big(A^q_n\big)_{ii},  
\eqno (1.3)
$$
and 
 $$ 
Z_n^{(Y^{(q)})}(g)= \log \bE_{ER} \exp \left\{- g Y_n^{(q)}\right\} 
, \quad Y_n^{(q)} = \sum_{i,j=1}^n \big( A^q_n\big)_{ij}
\eqno (1.4)
$$
have been studied in the limit $n\to\infty$, where $\bE_{ER}$ denotes the mathematical expectation with respect to the  measure generated by the family $\{\CA_n\}$ (1.1)
of Erd\H os-R\'enyi random graphs $\rG(n,p_n)$.
Using a diagram technique, it is shown in \cite{K-08} that in three asymptotic regimes of dense, dilute and sparse random graphs, in dependence of the rate of $p_n$, 
the  cumulants of properly normalized random variables $\hat X^{(q)}_n$ 
and $\hat Y^{(q)}_n$ converge,
$$
{1\over b_n} Cum_k(\hat X^{(q)}_n) \to \CK_k^{(q)}(\omega),
\quad {1\over d_n} Cum_k(\hat Y^{(q)}_n)\to \fL_k^{(q)}(\omega),
\quad n\to\infty, \  \omega = 1,2,3,
\eqno (1.5)
$$ 
with  coefficients $b_n$ and $d_n$ determined by $n$ and  $p_n$. Relations (1.5) can be 
useful in the proof of limit theorems either 
the  Central Limit Theorem or the Poisson one (1.2) \cite{K-08}. 


\vskip 0.2cm 
The classical Erd\H os-R\'enyi ensemble $G(n,c/n)$  can be regarded as the simplest model  of random graphs
because of the uniformity of the edge probability $p_n(i,j)=c/n$. 
In the present paper,
we consider an ensemble of modified Erd\H os-R\'enyi random graphs such that  
 edge probability $p_n(i,j)$ (1.1)  decreases when the difference $|i-j|$ does. 
These can be implemented by modifying  adjacency matrices $A_n$ 
by 
an additional factor,  
$$
(\hat A_n^{(c,R)})_{ij} = (A_{n}^{(c/n )})_{ij} w^{(R)}_{ij}, \quad 1\le i  < j\le n,
\eqno (1.6)
$$ 
where random variables $w^{(R)}_{ij}$ are independent from the family ${\CA_n}= \{ a_{ij}^{(n)}, 1\le i< j\le n\}$, 
$$
w^{(R)}_{ij} = \begin{cases}
 1,  &  \text{with probability $\sigma(|i-j|/R)$}  , \\
0,  \quad   & 
\text {with probability $1-\sigma(|i-j|/R)$
},
\end{cases}
\eqno (1.7)
$$
and the edge  probability is such that $\s(x)$ decreases to zero as $x\to\infty$.
Random graphs determined by adjacency matrices $\hat A_n^{(c,R)}$ (1.6)
 resemble certain  versions of  geometric random graphs (see \cite{NRS} and references therein). 
Also one can trace out a similarity  between  (1.6) and  random adjacency matrices of one-dimensional  long-range percolation radius model
\cite{BB,CGS}.
 
The interest of  studies  of large random graphs associated with $\{\hat A_n^{(c,R)}\}$ is that a new parameter $R$ appears
thus generalizing  the classical Erd\H os-R\'enyi ensemble.   
On could say that  $R$ represents a kind of effective raduis (or dimension)  of the graphs we consider.
In the limiting transition when  $n,c,R$ tend simultaneously to infinity, a new type of asymptotic behavior of variables 
$X$ (1.3) and $Y$ (1.4) can be observed. In particular, we have detected an asymptotic regime, when,
in contrast to (1.2),  the average vertex degree 
of the graph remains finite while the total number of triangles  tends to infinity.  
Also, one can find  an asymptotic regime, when the  number of triangles infinitely increases while the average number of quadrangles
tends to zero.
This cannot happen
in the Erd\H os-R\'enyi ensemble.

The paper is organized as follows. In Section 2, we  determine random  matrix ensemble  equivalent to (1.9), (1.10) 
and formulate our main results.
In Section 3, we develop a general diagram technique to study  cumulants of random variables that we consider. 
In Section 4, we prove  theorems 
for cumulants of the number of non-closed $q$-step walks $Y^{(q)}$. 
In Section 5, we prove theorems for cumulants of the number of closed $3$-step walks $X^{(3)}$. 
In Section 6,  
we show that either the  Central Limit Theorem
or the  Poisson Limit Theorem are valid for these  variables, 
in dependence of the asymptotic regime considered. 
These results allow us, in particular,  
to determine the asymptotic regime 
that leads to a rigorous  solution of  the random graph collapse problem.
In Section 7, we present  a  
color version of the Pr\"ufer codification procedure  adapted to the tree-type diagrams we consider and  deduce explicit expressions 
for limit  cumulants.
  In  Section 8, 
we 
describe analyticity properties of generating function of limiting cumulants of $Y$-models.


  \section{Main results}

Let us consider a family $\CA_{N,c,R}= \{a_{ij}^{(N,c,R)}, -n\le i< j\le n\}$ of jointly independent random variables
$$
a_{ij}^{(N,c,R)} = \begin{cases}
 1,  &  \text{with probability $p_{N,c,R}= {\ds c\over \ds N} e^{-\psi^2((i-j)/R)}$}, \\
0,  \quad   & 
\text {with probability $1- p_{N,c,R} 
$,
}
\end{cases}
\quad \ -n\le i< j\le n, 
\eqno (2.1)
$$
where 
$N= 2n+1$ and $\psi(x), x\ge 0$  is a continuous even real function.
Real symmetric $N$-dimensional  random matrices with elements
$(A_{N,c,R})_{ij} = a_{ij}^{(N,c,R)}$, $-n\le i< j\le n$ can  be considered
as the adjacency matrices of  random graphs $\g_N = (\rV_N, \rE_N)$, 
 where 
$\rV_N$ is the set of  $N$ ordered vertices  labeled by integers from 
$\rL_N= \{-n,\dots, n\}$ and $\rE_N$ 
is a  subset of  pairs $\{i,j\}$, $i,j \in \rV_N$, $i\neq j$. 
We denote by $\bE_{N,c,R}$ the mathematical expectation with respect to the measure generated by random variables 
(2.1).  

Let us introduce a real  symmetric  matrix 
with zero diagonal and the elements 
$$
\left( A^{(\alpha)}_{N,c,R}\right)_{ij} = 
\left(1 + \alpha \psi^2\left( { i-j\over  R}\right)
\right)\, a_{ij}^{(N,c,R)}, \quad -n\le i<j\le n,
\eqno (2.2)
$$
 We consider    random variables, 
$$
X^{(\a,q)}_{N,c,R}= \Tr \big(A^{(\a)}_{N,c,R}\big)^q,
\eqno (2.3)
$$
$$
Y^{(\a,q)}_{N,c,R}= \sum_{i,j\in \rL_N}
\left(\big(A^{(\a)}_{N,c,R}\big)^q\right)_{ij} 
\eqno (2.4)
$$
and study asymptotic behavior of their  cumulants
$$
\rcum_k \left(X^{(\a,q)}_{N,c,R}\right) = \rcum_k(X^{(q)}) 
\quad {\hbox{and}} \quad  
\rcum_k \left(Y^{(\a,q)}_{N,c,R}\right) = \rcum_k(Y^{(q)})
\eqno (2.5)
$$
in the limiting transition  
$$
N, c, R \to \infty, \quad R= o(N), \quad  c=o(N)
\eqno  (2.6)
$$ that we denote by $(N,c,R)_0\to\infty$. 

We start with random variables $ 
Y^{(\a,q)}_{N,c,R}$ (2.4).
It is suitable to formulate our results in the cases $\a=0$ and $\a=1$ of (2.2)
separately. 
Everywhere below, 
 we omit the subscripts and superscripts $N,c,R$ when no confusion can arise.

\vs 
{\bf Theorem 2.1.}  {\it  Let   $\psi(x), x\in \bR $  be an even continuous  strictly positive function that is monotone increasing for $x\ge 0 $
and such that 
$$
V_0= \int_{-\infty}^\infty e^{ -  \psi^2(t)} dt <\infty.
\eqno (2.7)
$$ 
Then for any given $k\in \bN$, the following limits exist: 

\vs 
\n i)  if $cR/N\gg 1$, then}
$$
\lim_{(N,c,R)_0\to\infty}\  {1 \over  c R} \rcum_k\left({N^{q-1}\over (c R)^{q-1}} Y^{(0,q)}\right)
= \Phi_k^{(q,1)}; 
\eqno (2.8) 
$$
{\it ii) if $cR/N=s>0$, then  
$$
\lim_{(N,R,c)_0\to \infty}\ {1\over cR} \, \rcum_k(Y^{(0,q)})= 
\Phi_k^{(q,2)}(s); 
\eqno (2.9)
$$}
{\it iii) if $ cR/N \ll 1$, then 
  }
$$
\lim _{(N,c,R)_0\to\infty}  {1\over cR} \rcum_k(Y^{(0,q)}) 
= \Phi_k^{(q,3)} = 2^{k-1} V_0, 
\eqno (2.10)
$$
{\it where} 
$$
\Phi_k^{(q,1)} = t_k^{(q)} V_0^{k(q-1)+1},
\quad 
t_k^{(q)}= 2^{k-1} \, q^k \, \big( k(q-1)+1\big)^{k-2},
\eqno (2.11)
$$
{\it and } 
$$
\Phi^{(q,2)}_k(s)= \sum_{l=1}^{k(q-1)+1} s^{l} \phi_k^{(q,l)},   
\eqno (2.12)
$$
{\it where}
$$
\phi_k^{(q,l)}= 
2^{k-1}\,  V_0^{l-k+1} \, 
\left( l - k+1\right)^{k-2} 
\sum_{\stackrel{1\le r_1, r_2, \dots, r_k\le q}
{r_1+\cdots+ r_k=l}} \ 
\prod_{i=1}^k \, r_i \rT^{(q)}(r_i).
\eqno (2.13)
$$

In (2.11),   $t_k^{(q)}$   represents the number of maximal 
tree-type diagrams constructed with the help of $k$ linear graphs with $q$ edges
that we denote by $\l_q$. 
We determine the tree-type diagrams in Section 4
and deduce explicit expressions for  $t_k^{(q)}$ 
in  Section 7.
In (2.13), 
$\phi_k^{(q,l)}$ 
represents the number 
of all possible tree-type diagrams constructed with the help of 
$k$ linear graphs with $r_i$ edges, $1\le r_i\le q$, 
with $r_1+\dots r_k=l$
(see Section 7 for the definition of  $\rT^{(q)}(r_i)$).
In (2.10),  factor $2^{k-1}$ represents the number of   minimal tree-type diagrams (see Section 4).
Summing up, we can say that the leading contribution 
to cumulants $\rcum_k(Y^{(q)})$ (2.5) is determined in  three asymptotic regimes 
(i), (ii), and (iii) 
by the maximal tree-type diagrams, all tree-type diagrams  and the minimal tree-type diagrams, respectively.
  


\vs 
{\bf Theorem 2.2.}  {\it Let $\psi(x), x\in \bR$ be as in Theorem 2.1. 
If
$$
V_{m}= \int_{-\infty}^\infty \big(1+ \psi^2(s)\big)^{m} \, e^{- \psi^2(s)}ds< \infty,
\quad \forall m\in {{\mathbb N}},
\eqno  (2.14)
$$
then the following limits exist for any $k\in \bN$:
\vs 
\n i) if  $cR/N\gg 1$, then }
$$
\lim_{(N,c, R)_0\to\infty}\  {1 \over c R }\,  
\rcum_k\left({N^{q-1}\over (c R)^{q-1}} Y^{(1,q)}\right) 
=  \Xi_k^{(q,1)} ;
\eqno (2.15)
$$

{\it ii) if $cR/N=s$, then  
$$
\lim_{(N,R,c)_0\to \infty}\ {1\over cR} \, \rcum_k\big(Y^{(1,q)}\big)= 
\Xi_k^{(q,2)}(s);
\eqno (2.16)
$$}

\vskip 0.2cm 
{\it iii) if $ c R/N \ll 1$, then }
$$
\lim _{(N,c,R)_0\to\infty}  {1\over cR} \rcum_k\big(Y^{(1,q)}\big) 
= \Xi^{(q,3)}_k=  2^{k-1} V_{kq}.
\eqno (2.17)
$$
{\it In (2.15),
$$
\Xi_k^{(q,1)} = 
{2^{k-1}\, q^k\,  (k-1)!\over k(q-1)+1} \ \sum_{u=1}^{k-1}  \, u! \,  {{k(q-1)+1}\choose{\kappa }}
\ \sum_{ \stackrel{ \s_k= (s_1, \dots, s_{k-1})} { | \s_k| = u, \, 
\Vert  \s_k\Vert = k-1}} 
\ V_1^{s_0}\ 
\prod_{i=1}^{k-1} {1\over s_i!} \left( {V_{i+1}\over i!}\right)^{s_i}, 
\eqno (2.18)
$$
 where  
 $$
 |\s_k |= s_1+\dots + s_{k-1},\ s_i\ge 0, \quad  
 \Vert \s_k\Vert = \sum_{i=1}^{k-1} i s_i,
 \quad  s_0= k(q-1)+1 - |\s_k|.
 \eqno (2.19)
 $$ 
}

{\it Remark 1.}
Results of Theorems 2.1 and 2.2 
presented in the third asymptotic regime 
are also valid in the case of constant concentration $c= Const$
and the limiting transition (2.6) replaced by  
$
N,R\to\infty$,  $R= o(N).
$


{\it Remark 2.} All statements of Theorem 2.1 and Theorem 2.2 remain true in the case when
 $\psi(x)=0$ and $R=N$;
 in this case, it is sufficient to replace  $V_{j}$ by 1  for all $j\ge 0$  (2.14). 
 


{\it Remark 3.} Relation (2.18) will be proved in Section 7. 
The right-hand side of  (2.16) can be  computed explicitly
with the help of  
 (2.12), (2.13) and (2.18).
General expression is  cumbersome;  we present its value for   the first and the 
 second cumulants
 in the case of  $q=2$,
$$
\Xi_1^{(2,2)}(s)= sV_1+V_2, \quad 
\Xi_2^{(2,2)}(s)= 
8s^2 V_1^2V_2   + 8 s V_1V_3  + 2 V_4
\eqno (2.20)
$$
(see Section 7 for the details).

Let us pass to 
the number  of closed walks 
 $X^{(\a,q)}$ (2.3).
Our primary interest is related with  the number of triangles  $T_n$ (1.2) 
and we consider
$X^{(\a, q)}$-models with $q=3$ only.

{\bf Theorem 2.3.} 
 {\it Under conditions of Theorems 2.1 and 2.2, 
 there exist numbers 
 $\Theta_k^{(\a,\ri)}$,  $\Theta_k^{(\a,\rii)}$
and 
 $\Theta_k^{(\a,\riii)}$ such that 
\vskip 0.1cm  
i) if 
$c^2R/N^2 \gg 1$, then  
}
$$
\lim_{(N,c, R)_0 \to \infty} {1\over   c R} \rcum_k
\left( {N^2\over c^2 R}  X^{(\a,3)}\right)= 
\Th_k^{(\a,\ri)};
\eqno  (2.21)
$$

 {\it ii) if 
$c^2R/N^2 = s$, then 
 } 
$$
\lim_{(N,c, R)_0 \to \infty} {1\over   c R} \rcum_k
\big(  X^{(\a,3)}\big)= 
\Th_k^{(\a,\rii)}(s);
\eqno (2.22)
$$

{\it iii) $c^2R/N^2\ll 1$, then 
$$
\lim_{(N,c, R)_0 \to \infty} {N^2\over c^3 R^2} 
\rcum_k\big( X^{(\a, 3)}\big)= 
\Th_k^{(\a,\riii)}=
 {6^{k-1}} \tilde H^{(\a,3)}_k, 
 \eqno (2.23)
$$ 
where we denoted
$$
\tilde H^{(\a,3)}_k = 
{1\over 2\pi} \int_{-\infty}^\infty \big(\tilde  h_k^{(\a)}(p)\big)^3 dp
\eqno (2.24)
$$
and
$$
\tilde h_k^{(\a)}(p) = \int_{-\infty} ^\infty h_k^{(\a)}(x) e^{-ipx} dx, 
\quad 
h_k^{(\a)}(x)= (1 + \alpha \psi^2(x))^k e^{ - \psi^2(x)}.
$$
}

{\it Remark 1.} In analogy with Theorems 2.1 and 2.2,
terms $\Theta^{(\a,\ri)}_k$, $\Theta^{(\a,\rii)}_k(s)$ and 
$\Theta^{(\a,\riii)}_k$ represent total contributions of maximal tree-type diagrams, 
all tree-type diagrams and minimal tree-type diagrams to the cumulants of 
$X^{(\a,3)}$-model (2.5), respectively.
In  Section 5 we give rigorous definition of tree-type diagrams
for $X^{(\a,3)}$-models.

{\it Remark 2.} Expressions for $\Th^{(\a,\ri)}_k$ and $\Th^{(\a, \rii)}_k(s)$ with 
both $\a=0$ and $\a=1$ 
are rather complicated
and involve \mbox{products} of  functions $\tilde h_m^{(1)}(p)$
and their convolutions. 
We present limiting expressions 
for the first cumulant and the second cumulant only; these are
$$
\Th_1^{(\a,\ri)}= {1\over 2\pi} \int_{-\infty}^\infty
   \big(\tilde h_2^{(\a)}(p)\big)^3 dp, \quad 
\Th_2^{(\a, \ri)}= {9\over \pi}\int_{-\infty}^\infty 
 \left( 
\big( \tilde h_1^{(\a)}\big)^2*\tilde h_2^{(\a)}\right)(p)\  
\big( \tilde h_1^{(\a)}(p)\big)^2 dp , 
\eqno (2.25)
$$
and 
$$
 \Th_1^{(\a,\rii)}(s)= s \Th_1^{(\a,\ri)}, \quad 
 \Th_2^{(\a,\rii)}(s)= s^2 \Th_2^{(\a,\ri)} + 
 3s \Th_1^{(\a,\ri)}.
 \eqno (2.26)
$$


{\it Remark 3.} All statements of Theorem 2.3 remain true 
in the case when   $\psi^2(x)=0$ and   $R=N$; in this case, it is sufficient to replace $V_j$ by 1 for all $j\ge 0$. 
Then the  right-hand side of (2.21) 
is expressed by the number of tree-type diagrams (2.11) with $q=3$, 
$$
\Th_k^{(\ri)}= t_k^{(3)}= 2^{k-1}3^{k} (2k+1)^{k-2},\quad k\ge 1
\eqno (2.27)
$$
and
the right-hand side of (2.22) is given by 
$$
\Th_k^{(\rii)}(s)= 
{k!\over 3} \sum_{u=1}^k   
\, s^u\,  6^{u-1} (2u+1)^{u-2}
\sum_{ \stackrel{ \s_{k+1}= (s_1, \dots, s_{k})} { | \s_{k+1}| = u, \, 
\Vert  \s_{k+1}\Vert = k}}  \prod_{i=1}^{k} {1\over (i!)^{s_i} s_i!}, 
\eqno  (2.28)
$$
where the sum runs over variables $\s_{k+1}$ of the form (2.19). 
 
   
\section{Diagrams and their contributions}

Asymptotic behavior of cumulants (2.5)
 can be studied with the help of  well-known technique of connected diagrams. 
 In this section, 
we consider  random variables $Y^{(\a,q)}_{N,c,R}$ (2.4), 
that can be presented  as follows,
$$
Y^{(\a,q)}_{N,c,R}= \sum_{i_1=-n}^n \dots \sum_{i_{q+1}=-n}^n 
\CY^{(\a, q)}\big(\la  i\ra _{q+1}\big), \quad \CY^{(\a, q)}\big(\la  i\ra _{q+1}\big)
=a^{(\a)}_{i_1 i_2} a^{(\a)}_{i_2 i_3} \cdots a^{(\a)}_{i_q i_{q+1}},
\eqno (3.1)
$$
where 
$
a^{(\a)}_{ij}= (A^{(\a)}_N)_{ij},  
$
and where we denoted 
$\la  i\ra_{q+1} = (i_1, \dots , i_{q+1})$.
In what follows, we omit the superscripts $\a$ in $Y^{(\a,q)}$ and $\CY^{(\a,q)}$ 
when no confusion can arise.  
After obvious modifications,  arguments of this section can be applied 
 to random variables
$X^{(\a,q)}_{N,c,R}$ (2.3).
  
\subsection{Cumulants and mixed cumulants }
 We start with  a  representation of cumulants $\rcum_k(Y^{(q)})$, 
  $Y^{(q)}=Y^{(q)}_{N,c,R}$
 by mixed cumulants  $cum_k\left\{ Y^{(q)}, \dots, Y^{(q)}\right\}$,
 given by equality 
 $\rcum_k(Y^{(q)}) = cum_k\left\{ Y^{(q)}, \dots, Y^{(q)}\right\}$ (see  \cite{MM} and more recent monograph \cite{PT}), 
  where by definition
  $$
  cum_k\left\{ \Theta_1, \dots, 
\Theta_k\right\}= {d^k\over dz_1 \cdots dz_k} 
\log 
\bE \left(  \exp\left\{z_1 \Theta_1+ \dots+ z_k\Theta_k
\right\}\right) \vert_{z_1=0,\dots, z_k=0}.
\eqno (3.2)
$$
Then $$\rcum_k(Y^{(q)})=
\sum_{
\big( \la i \ra^{(1)}_{q+1},\,  \dots \, ,  \la i\ra^{(k)}_{q+1}\big)} 
cum_k\left\{ 
\CY^{(q)}\big(\la i \ra^{(1)}_{q+1}\big), \dots,  \CY^{(q)}\big(\la i\ra^{(k)}_{q+1}\big)
\right\},
\eqno (3.3)
$$
where we have used 
the multi-linearity property of 
mixed cumulants 
$$
cum_k\left\{a'\Th' + a''\Th'', \Th_2, \dots,\Th_k\right\}=
a' cum_k\left\{\Th', \Th_2, \dots,\Th_k\right\}
+a''cum_k\left\{\Th'', \Th_2, \dots,\Th_k\right\}
$$
that is an easy exercice based on definition (3.2)
(see \cite{MM}, Chapter 2).
In the right-hand side of (3.3), it is useful to consider 
$\la i\ra ^{(j)}_{q+1}$ as a realization 
given by  an element of the vector space $\rL_N^{q+1}$, i.e. 
as  
 $i^{(j)}_l, 1\le l\le  q+1$ with given values of variables
that belong to the set $\rL_N=\{-n,\dots, n\}$. We denote this 
realization by  
$$
\la i\ra^{(j)}_{q+1} = \la i_1^{(j)}, \dots, i^{(j)}_{q+1}\ra_N,
\quad -n \le i^{(j)}_l\le n, \quad l=1, 2, \dots, q+1
\eqno (3.4)
$$
and write $\la \CY^{(q)}\ra_j$ instead of $\CY^{(q)}\big(\la i\ra^{( j)}_{q+1}\big)$ and  consider the sum of (3.3) as the one running through 
 the set  of all possible such realizations 
 $\rL_N^{k(q+1)}$.

The right-hand side of (3.3) can be computed
with the help of the following
 \mbox{formula 
\cite{MM,PT},}
$$
cum_k\left\{\la \CY^{(q)}\ra_1, \dots, \la\CY^{(q)}\ra_k\right\} = 
\bE \left(\la \CY^{(q)}\ra_1  \cdots 
\la\CY^{(q)}\ra_k \right)
$$
$$
+ \  
\sum_{s=2}^k \ 
\sum_{\pi_s \in \Pi_k}\,  (-1)^{s-1}\,  (s-1)! \ 
\bE \left(\la \BY^{(q)}(\rho_1)\ra\right)\,  \cdots 
\bE\,  
\left( \la \BY^{(q)}(\rho_s)\ra\right),
\eqno (3.5)
$$ 
where $\pi_s = (\rho_1, \dots, \rho_s)$ is a partition of the set 
$\{1, 2, \dots, k\}$ into $s$ subsets, $\Pi_k$ is a family of all possible partitions and 
$$
\bE \la \BY^{(q)}({\rho_l}) \ra= 
\bE \prod_{j\in \rho_l} \la \CY^{(q)}\ra_j.
\eqno (3.6)
$$
The subsets $\{\rho_l\}_{1\le l\le s}$ can be ordered by recurrence, the minimal element of $\rho_{i+1} $ being the minimal one
that is not contained in $\rho_1\cup \dots\cup \rho_i$. 

To compute the  mean values of the form  (3.6), we use a version of the diagram method 
widely known in random matrix theory and its applications. 
Some of its elements 
date back to the pioneering  works by 
E. Wigner \cite{W}.
The starting point   is to use a 
graphical representation of   random variable $a_{ij}$ by 
two vertices joined by an edge, the vertices being attributed by 
values of $i$ and $j$, respectively. Given a realization
(3.4),  
the product 
$
\la \CY^{(q)}\ra= a^{(\a)}_{i_1 i_2} a^{(\a)}_{i_2 i_3} \cdots 
a^{(\a)}_{i_q i_{q+1}} 
$
 can be  represented by an ensemble of vertices and edges 
that we denote by $\rg(\la i\ra_{q+1})$. 
If the sequence $\la i\ra_{q+1}= \la i_1, i_2, \dots, i_{q+1}\ra_N$ has two or more elements with equal value,
then we draw one vertex $v$ attributed by this value. This means that 
 $\rg(\la i\ra_{q+1})$ 
is given by  a diagram of the form of a (multi)-graph 
with  $q$ edges when counted with its multiplicities. 
The key observation of this technique is that different multi-edges of $\rg(\la i\ra_{q+1})$
represent independent random variables and then 
$$
\bE\la  \CY^{(q)}\ra
= 
\prod_{e(j,l)\in \CE(\rg(\la i\ra_{q+1}))}
\bE \left(\left( a^{(\a)}_{i_ji_l}\right)^{m(j,l)}\right),
\eqno (3.7)
$$
where  $\CE (\rg(\la i\ra_{q+1}))$ is  the ensemble of edges of 
$\rg(\la i\ra_{q+1})$ and $m(j,l)$ is the multiplicity of the edge 
$e(j,l)= (i_j, i_l)$ of $\CE$. 
By these rules,  the product 
$\la \CY^{(q)}\ra_1 \la \CY^{(q)}\ra_2 \dots \la \CY^{(q)}\ra_k$
can also be  represented by a multi-graph  
and its average value can be computed with the help of (3.7).
In the next subsection, we give rigorous definition of the 
diagrams that determine the classes of equivalence we need. We follow mostly the lines of \cite{K-08}.



\begin{figure}[htbp]
\centerline{\includegraphics
[height=16cm, width=25cm]{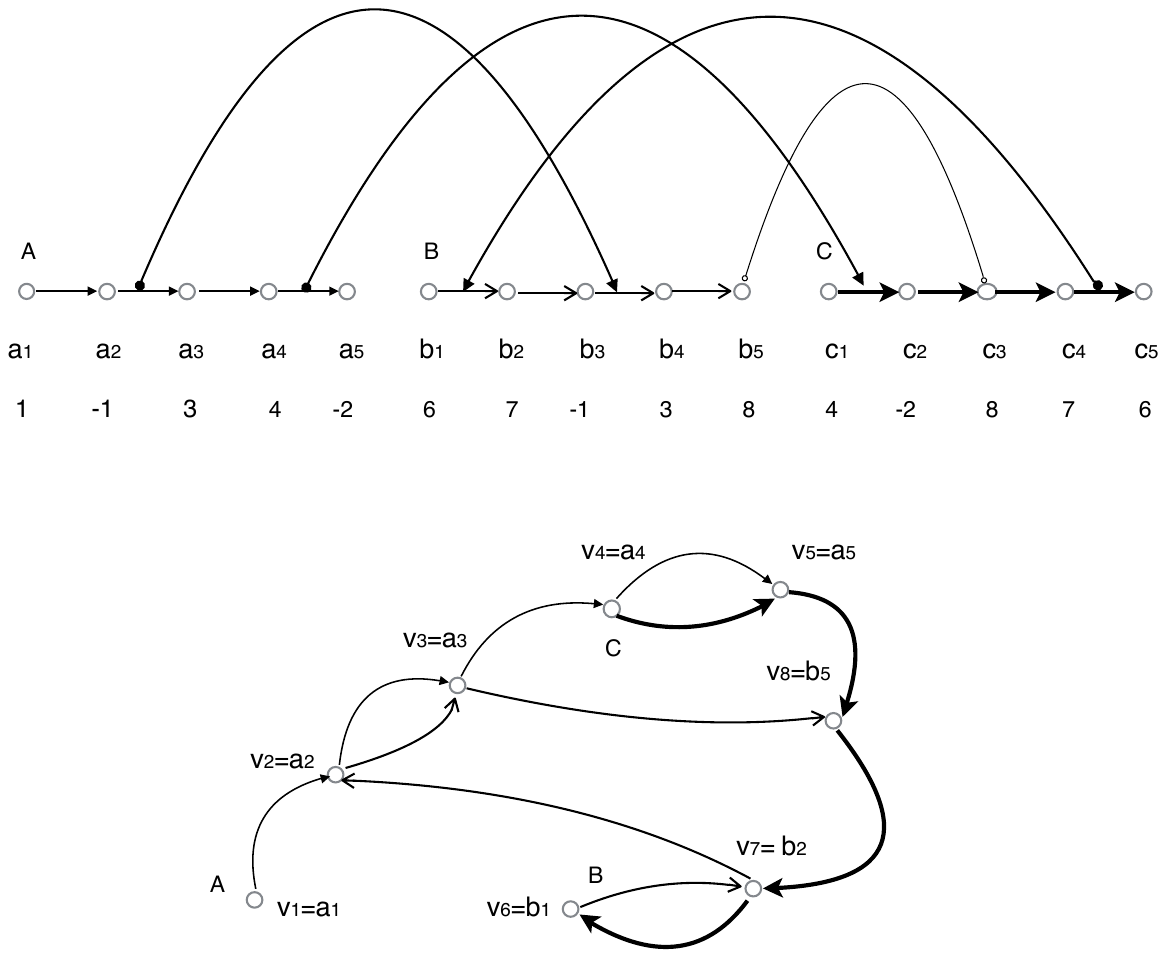}}
\caption{{{\it Type-I-diagram $\CD_3^{(4)}$ and type-II-diagram $D_3^{(4)}$ of the same $\CL_3^{(5)}(N)$}}}
\end{figure}


\subsection{Diagrams of $\l_q$-elements} 

Regarding the case of $Y$-models, we describe the diagram technique to compute the mathematical expectation
$\bE \left(\la \CY^{(q)}\ra_1  \cdots 
\la\CY^{(q)}\ra_k \right)$ of (3.5) and 
 start by  drawing  $k$ linear graphs 
$\l_q^{(j)}, j=1, \dots, k$, each $\l_q$ having  $q+1$ vertices and $q$ oriented in the same direction edges. We refer to $\l_q^{(j)}$ as to the $\l$-elements and  denote the set of ordered $\l$-elements by $\L_k^{(q)}= \{ \l_j^{(q)}, 1\le j\le k\}$. 


We denote vertices of  $\l_q^{(j)}$ by
$\u_1^{(j)}, \u_2^{(j)}, \dots, \u_{q+1}^{(j)}$ and denote the edges of 
$\l_q^{(j)}$ by $(\u^{(j)}_l, \u^{(j)}_{l+1})$. 
Regarding   a realization
 $$
\CL_k^{(q+1)}(N)= \left( \la i\ra^{(1)}_{q+1},\ \la i\ra^{(2)}_{q+1}, \   \dots \  , 
\la i\ra^{(k)}_{q+1}\right)_N, 
\eqno (3.8)
$$
we attribute to each vertex $\u_l^{(j)}$ the number 
$\la i_l^{(j)}\ra \in \rL_N$. 
Thus one can speak about  a realization of $\l$-elements 
$ \la \l_q^{(1)}, \dots, \l_q^{(k)}\ra_N$ of the set of ordered elements and   
denote by $\la \u_l^{(j)}\ra $  a realization of a  vertex
i. e.   the vertex together with the number attributed to it.
We denote this realization by $\la \L_k^{(q)}\ra_N= (\L_k^{(q)}, \CL_k^{(q+1)}(N))$. 

If realization $\la \l_q^{(1)}, \dots, \l_q^{(k)}\ra_N$ is such that
there exist  edges $e= (\u, \theta)$ and $e'=(\u',\theta')$  
that 
 either  $\la \u\ra = \la \u'\ra $, $\la \theta\ra =\la \theta'\ra $ or
 $\la \u\ra = \la \theta'\ra $, $\la \theta\ra = \la \u'\ra $, then  
 we say that these  edges are positively or negatively coupled. 
Regarding a realization $\CL_k^{(q+1)}$, we construct  a diagram of the first type $\CD_k^{(q)}= \CD(\CL_k^{(q+1)})$
and a diagram of the second type $D_k^{(q)}= D(\CL_k^{(q+1)})$, whose rigorous definitions will be presented below. 

 Type-I-diagram $\CD_k^{(q)}= \CD(\CL_k^{(q+1)})$ is given by  set $\L_k^{(q)}$ of $\l$-elements, where 
 positively and negatively coupled edges of $\la \L_k^{(q)}\ra_N$ are joined by positively or negatively oriented arcs;  a positively oriented arc goes from the edge of the minimal $\l$-element of the corresponding pair of $\l$-elements to the maximal $\l$-element of this pair; a  negatively oriented arc goes in  opposite direction. We refer to these arcs as to e-arcs. The diagram $\CD_k^{(q)}$ contains also non-oriented arcs that we refer to as v-arcs; these arcs
 the vertices attributed by the same number from $\rL_N$ but not attached to the edges joined by e-arcs.
 The collection of e-arcs is minimized by 
keeping the arcs between the edges of nearest 
 of $\l$-elements and erasing all other e-arcs; the collection of v-arcs is minimized by the same procedure.
  Let us note that the diagram $\CD_k^{(q)}$ obtained is constructed with the help of realization $\CL_k^{(q+1)}$ (3.8) but does not contain it, there can be several realizations leading to the same diagram $\CD_k^{(q)}$.



Type-II-diagram $D_k^{(q)}= D(\CL_k^{(q+1)})$ is given by  set $\L_k^{(q)}$ of $\l$-elements, where 
vertices attributed by the same numbers in $\la \L_k^{(q)}\ra_N$ are identified. 
Therefore $D_k^{(q)}$ can be determined as a multi-graph obtained from type-I-diagram  $\CD_k^{(q)}$ by gluing edges and vertices.

We denote by $\CV(D_k^{(q)})$ and 
 $\CE (D_k^{(q)})$ the sets of  vertices and  edges of $D_k^{(q)}$, respectively
 and denote by 
 $
 V= | \CV(D_k^{(q)})|$ 
the cardinality  of $\CV(D_k^{(q)})$. Given  an edge $e = (\rv', \rv'') \in \CE(D_k^{(q)})$, we denote by 
$m(e)$ its multiplicity.
The vertices of $\CD_k^{(q)}$ are labeled 
by letters v$_i$ in  the order prescribed by the order of  vertices of $ \CD(\CL_k^{(q+1)}(N))$.
We do this by recurrence: we set v$_1=\u_1^{(1)}$ and v$_2=\u^{(1)}_2$; having $l$ vertices v$_1$, \dots, v$_l$ labeled, 
we attribute the label v$_{l+1}$ for the next in term vertex of $D(\CD_k^{(q)})$, according to the order of appearance in $\CD_k^{(q)}$, that  
is not already labeled. 
Having all vertices labeled, we get  the type-II-diagram 
$D(\CL_k^{(q+1)}(N))$.

Let us note that  $\l$-elements of   type-II-diagram keep their ordering and  orientation of   multi-edge components; so, 
there is a bijection between 
type-II-diagram $\CD(\CL_k^{(q+1)})$ and type-I-diagram $\CD(\CL_k^{(q+1)})$; these diagrams of two types  give different graphical 
representation 
of the same realization $\CL_k^{(q+1)}$ and are, in fact, equivalent but useful in different situations. 

Regarding a  type-II-diagram $D_k^{(q)}$ and replacing each multi-edge by a non-oriented 
simple edge, we get a graph $G_k^{(q)}= G(D_k^{(q)})$.   It is clear that $\CV(G_k^{(q)})= \CV(D_k^{(q)})$. 

On Figure 1, we present a realization 
type-I-diagram
$\CD_3^{(4)}= \CD(\CL_3^{(5)}(N))$ and type-II-diagram $D_3^{(4)}= D(\CL_3^{(5)}(N))$
determined by 
$$
\CL_3^{(5)}(N)= \left( \la 1,-1,3,4,-2\ra_N, \la 6,7,-1,3,8\ra_N, \la 4,-2,8,7,6\ra_N\right);
$$ 
for simplicity,  the elements $\l_4^{(1)}, \l_4^{(2)}, \l_4^{(3)}$ are marked by letters $\rmA,\rmB,\rmC$, respectively, and 
the vertices of $\l$-elements are labeled   by $\rma_i,\rmb_j,\rmc_l$, $1\le i,j,l\le 5$. 



We say that two realizations $
 \CL_k^{(q+1)}(N)$ and $
\tilde  \CL_k^{(q+1)}(N)$ belong to the same class of equivalence, 
if their  type-II-diagrams coincide, 
$D(\CL_k^{(q+1)}(N)) =  D(\tilde  \CL_k^{(q+1)}(N))$.
Given a diagram $D_k^{(q)}$, we denote corresponding class of equivalence by $\CC(D_k^{(q)})$
and determine 
the weight of this diagram $W_{N,c,R}^{(\a)}(D_k^{(q)})$ as follows,
$$
\rW_{N,c,R}^{(\a,q)}(D_k^{(q)})=
\sum_{
\CL_k^{(q+1)}(N)
\in\,  \CC(D_k^{(q)}) 
} \bE 
 \left(
 \la \CY^{(q)}\ra_1
 \  \cdots \ 
\la \CY^{(q)}\ra_k
 \right).
 \eqno (3.9)
$$
The cardinality of the class of equivalence of $D_k^{(q)}$ is given by 
$N(N-1)\dots (N-V+1)$, where $V= |\CV(D_k^{(q)}|$. 
Taking into account  (3.7), we get from  (3.9)  the following equality for 
the weight of $D_k^{(q)}$ for a general value of $\a$,
$$
\rW_{N,c,R}^{(\a,q)}
(D_k^{(q)})=\sum_{\la s_1, \dots, s_V\ra^*_N}\ 
\prod_{\stackrel {1\le i<j\le V:}
{ e_{ij}= (\rv_i,\rv_j)\in \CE(D_k^{(q)})}
}\ 
{c\over N}\, \,  h^{(\a)}_{m(i,j)}\left( {s_i-s_j\over R}\right), 
\eqno (3.10)
$$
where the sum is performed over all possible realizations of the $V$-plets $\la s_1, \dots, s_V\ra^*_N$ such that $s_i\neq s_j$ for all $(i,j)$ and 
where,  according to (2.1) and (2.2), 
$$
 h_{m(i,j)}^{(\a)} \big(t\big) = 
\big( 1+\alpha \psi^2(t)\big)^{m(i,j)} e^{- \psi^2(t)}, \quad m(i,j)=m(e_{ij}).
 $$
In what follows, we omit the superscripts 
$\a$ and $q$ in $\rW^{(\a,q)}_{N,c,R}$ when no confusion can arise. 
We see that the sum over all possible realizations (3.8) can be transformed
into the sum over  diagrams and we can write the following equality that involves the first term 
of the right-hand side of (3.5), 
$$
\sum_{\CL_k^{(q+1)}(N) \in \rL^{k(q+1)}_N
} 
\bE \left(\la \CY^{(q)}\ra_1\  \cdots 
\la \CY^{(q)}\ra_k\  \right)
=
\sum_{D_k^{(q)} \in \, \fD_k^{(q)}} \rW_{N,R}(D_k^{(q)}),
$$
where $\fD_k^{(q)}$ is the family of all possible type-II-diagrams. 
The last sum can be considered as the sum over diagrams $\CD_k^{(q)}$  belonging to the set of all 
type-I-diagrams 
 that we denote by the same letter $\fD_k^{(q)}$ as for the set of all possible type-II-diagrams. 



Regarding the right-hand side of (3.5), we observe that  $\pi_1$ represents the  trivial partition of the set $\{1, 2, \dots, k\}$ that consists of one subset only.
Then we can   write that 
$\rW_{N,c,R}(D_k^{(q)})= \rW_{N,c,R}(\CD_k^{(q)})=\rW_{N,c,R}(\CD_k^{(q)}(\pi_1))$. 
To study  terms of the right-hand side of (3.5)
with non-trivial partitions $\pi_s, s\ge 2$, we transform   $\CD_k^{(q)}$  into  a diagram 
$\CD_k^{(q)}( \pi_s)$;
 if an e-arc of $\CD_k^{(q)}$ joins  $\l$-elements that belong to different subsets $\rho$ and $\rho'$ of $\pi_s$, 
 then we depict this arc by  a dotted one; if an e-arc joins the $\l$-elements 
that belong to the  same subset $\rho$ of $\pi_s$,  then we draw it by an unbroken arc as in the case of trivial partition $\pi_1$. 
Remembering (3.2) 
  and  (3.5), we can write 
  that 
  $$
\rcum_k(Y^{(q)}) = \sum_{\CD_k^{(q)}\in \, \fD_k^{(q)}} 
\left( \rW_{N,c,R}\big(\CD_k^{(q)}\big) +   
\sum_{s=2}^k (-1)^{s-1} (s-1)!  \sum_{\pi_s \in \, \Pi_k} 
\rW_{N,c,R} \big( \CD_k^{(q)}(\pi_s)\big)\right) ,
\eqno (3.11)
$$
where, according to (3.6) and (3.9),
$$
\rW_{N,c,R} 
\big(\CD_k^{(q)}(\pi_s)\big) = \ 
\sum_{\CL_k^{(q+1)}(N)
\,  
\in\,  \CC(\CD_k^{(q)}(\pi_s)) } \ \, 
\prod_{l=1}^s
\  \bE \la \BY^{(q)}(\rho_l)\ra, 
\quad \BY^{(q)}(\rho_l)= 
\prod_{i_j\in \rho _l} \la {  \CY}^{(q)}\ra_{i_j} ,
\eqno  (3.12)
$$
and  subsets $\rho_l\subset \{1, 2, \dots, k\}$ are determined
by partition $\pi_s$. 

\subsection{Connected diagrams, weights and contributions}

We say that  diagram $\CD_k^{(q)}$ is  connected, if there is no  subset 
$\rho$ of $\{1, 2, \dots, k\}$ such that there is no e-arc that joins a $\l$-element of $\rho $ with
a $\l$-element of $\{1, 2, \dots, k\}\setminus \rho $. 
We denote connected type-I-diagrams 
by $\hat \CD_k^{(q)}$. 
 Non-connected 
type-I-diagrams 
will be denoted by   $\ddot  \CD_k^{(q)}$.

\vskip 0.2cm 

{\bf Lemma 3.1.} {\it For any non-connected diagram $\ddot  \CD_k^{(q)}$ the following relation is true, 
$$
\rW_{N,c,R}\big(\ddot  \CD_k^{(q)}\big)+ \sum_{s=2}^k \, (-1)^{s-1}\,  (s-1)! \sum_{\pi_s \in \Pi_k}  
\rW_{N,c,R}\big(\ddot   \CD_k^{(q)}(\pi_s)\big)
 = 0.
\eqno  (3.13)
$$
} 
{\it Proof.} 
If $\ddot  D_k^{(q)}$ is non-connected, then there exist at least two subsets 
$\rho'$, $\rho''$ 
 such that 
 $
 \rho'\cup \rho'' = \{1, \dots, k\}$
 and
   $\rho'\cap \rho''= \emptyset
  $
and such that  any element $\l_i, i\in \rho'$ is not connected to an element 
$\l_j, j\in \rho''$.  Therefore for each  realization 
$\CL_k^{(q)}(n)$ such that 
$\CD\big(\CL_k^{(q)}(n)\big)= \ddot \CD_k^{(q)}$  random variables 
$
\la \BY^{(q)}({\rho'})\ra  $ and 
$\la \BY^{(q)}({\rho''})\ra $ are jointly independent. Regarding the right-hand side of (3.3), we  can write  that
$$
\bE \left( \exp\big(z_1 \la \CY^{(q)}\ra_1+ 
\dots+ z_k\la \CY^{(q)}\ra_k\big)\right)
= \bE \left(\prod_{i\in \rho'} 
\exp \big(z_i \la \CY^{(q)}\ra_i\big)\right)
\, \bE \left( \prod_{i\in \rho''} 
\exp \big(z_i \la \CY^{(q)}\ra_i\big)\right)
$$
and thus   $cum_k(\la \CY^{(q)}\ra_1, \dots, \la \CY^{(q)}\ra_k)=0$.
Then  
$$
\sum_{\CL_k^{(q)}(N)\in \CC (\ddot \CD_k^{(q)})} 
cum_k(\la \CY^{(q)}\ra_1, \dots, \la \CY^{(q)}\ra_k)=0
$$
and  (3.13) follows. 
 $\Box$

Lemma 3.1 says that we can  restrict our consideration of the sum in  the right-hand side of (3.5)
to the ensemble 
of connected diagrams $\hat \CD_k^{(q)}$.
Let us first consider the case of the trivial partition $\pi_1$. In this case,  we can use 
 equivalence $\hat \CD_k^{(q)}(\pi_1)= \hat \CD_k^{(q)}= \hat D^{(q)}_k= \hat D_k^{(q)}(\pi_1)$. 

\vskip 0.2cm 
{\bf  Lemma 3.2.} 
{\it The order of the weight of diagram $\hat D_k^{(q) }(\pi_1)$
is given by the following equality, }
$$
\rW_{N,c,R} (\hat D_k^{(q)}(\pi_1))=
 O\left({c^{E } R^{V -1}\over N^{E-1}}\right), \quad (N,c, R)_0\to\infty,
\eqno  (3.14)
$$
{\it where   $E= | \CE(\hat D_k^{(q)})|$ and 
$V =
|\CV(\hat D_k^{(q)})|$.} 

{\it Proof.} 
We proceed  by recurrence. Let $\rv_{\max} $ be  the maximal 
 vertex of $\hat D_k^{(q)}$.
 Assume that there are $l$ edges $(\rv_{j_i},\rv_{\max})$, $1\le i\le l$ 
 attached to $\rv_{\max}=\rv_V$. We denote by $m(j_i,V)$ the multiplicity 
 of the edge $(\rv_{j_i}, \rv_{\max})$ and attribute variables $s_V$ and 
 $s_1, \dots, s_l$ to the vertices $v_V$ and $v_{j_1}, \dots, 
 v_{j_l}$, respectively.
 According to (3.10),
  the weight $\rW_{N,c,R}(D_k^{(q)})$
contains the following factor 
$$
P(s_1, \dots, s_l)= \prod_{i=1}^l \ \sum_{s_V=-n}^n \  
{c\over N} \ h_{m(j_i,V)}^{(\a)}
\left( {s_i-s_V\over R}\right). 
\eqno (3.15)
$$
Taking into account  the upper bound
$
\RM_{j_i} = \sup_{x\in \bR} \left( 1 + \psi^2(x)\right)^{m(j_i,V)} e^{-\psi^2(x)},
$
 we can write that 
$$
\sup_{s_2, \dots , s_l}  R^{-1} P(s_1, s_2, \dots, s_l)
\le \left({c\over N}\right)^l 
\ \CH_{m(j_1, V)}^{(R,\a)}\  \prod_{i=2}^l \RM_{j_i},
\eqno (3.16)
$$
where 
$$
\CH_{m(j_1, V)}^{(R,\a)} =
{1\over R} \, \sum_{s\in {\bZ}} h^{(\a)}_{m(j_1,V)}\left( {s\over R}\right).
$$
Using elementary inequalities
$$
\CH_m^{(R, \a)}\le \CH_m^{(R,1)} \le {1\over R} \sum_{i=0,1} (1+\psi^2(x_i))^m e^{-\psi^2(x_i)}
+ 2 \int_1^\infty (1+\psi^2(x))^m e^{-\psi^2(x)} dx,
$$
where $x_0=0$ and $x_1\ge 0$ is such that 
$$
(1+\psi^2(x_1))^m e^{-\psi^2(x_1)} = \sup_{x\ge 1} \left( (1+\psi^2(x))^m e^{-\psi^2(x)}\right),
\eqno (3.17)
$$
we conclude that erasing from $\hat D_k^{(q)}$ the vertex $v_{\max}$ together with 
all multi-edges attached to it produces a new diagram 
$\hat D_{k-\chi-m(j_1,V)}^{(q)}$,
with $\chi= m(j_2,V)+ \dots + m(j_l,V)$  whose weight is multiplied by 
$p_N^lR= (c/N)^l R$ 
and a constant bounded above by 
$\RM_\chi  \bar \CH_{m(j_1,V)}$,
where 
$\bar \CH_{m(j_1,V)} = \sup_R \CH_{m(j_1, V)}^{(R)}$. 
Here we have used an elementary  bound,
$
\RM_{j_2}\cdots \RM_{j_l}\le \RM_\chi
$. 
Regarding $\hat D_{k-\chi-m(j_1,V)}^{(q)}$, we repeat  the procedure described above. On the last step of this recurrence, 
we get the vertex $\rv_2$ attached 
to $\rv_1$ by a multi-edge $e(\rv_1,\rv_2)$. It is clear that this 
diagram with two vertices has a weight 
bounded by $NR (c/N)$ multiplied by $\bar \CH_{m(1,2)}$. 
Then we conclude that the total weight of the diagram 
$\hat D_k^{(q)}$ is of the order $(c/N)^{E} NR^{V-1}$,
where $E$ is the sum of all values of $l$ considered on each step of recurrence. 
This observation completes the proof of Lemma 3.2. $\Box$


\begin{figure}[htbp]
\centerline{\includegraphics
[height=19cm, width=15cm]{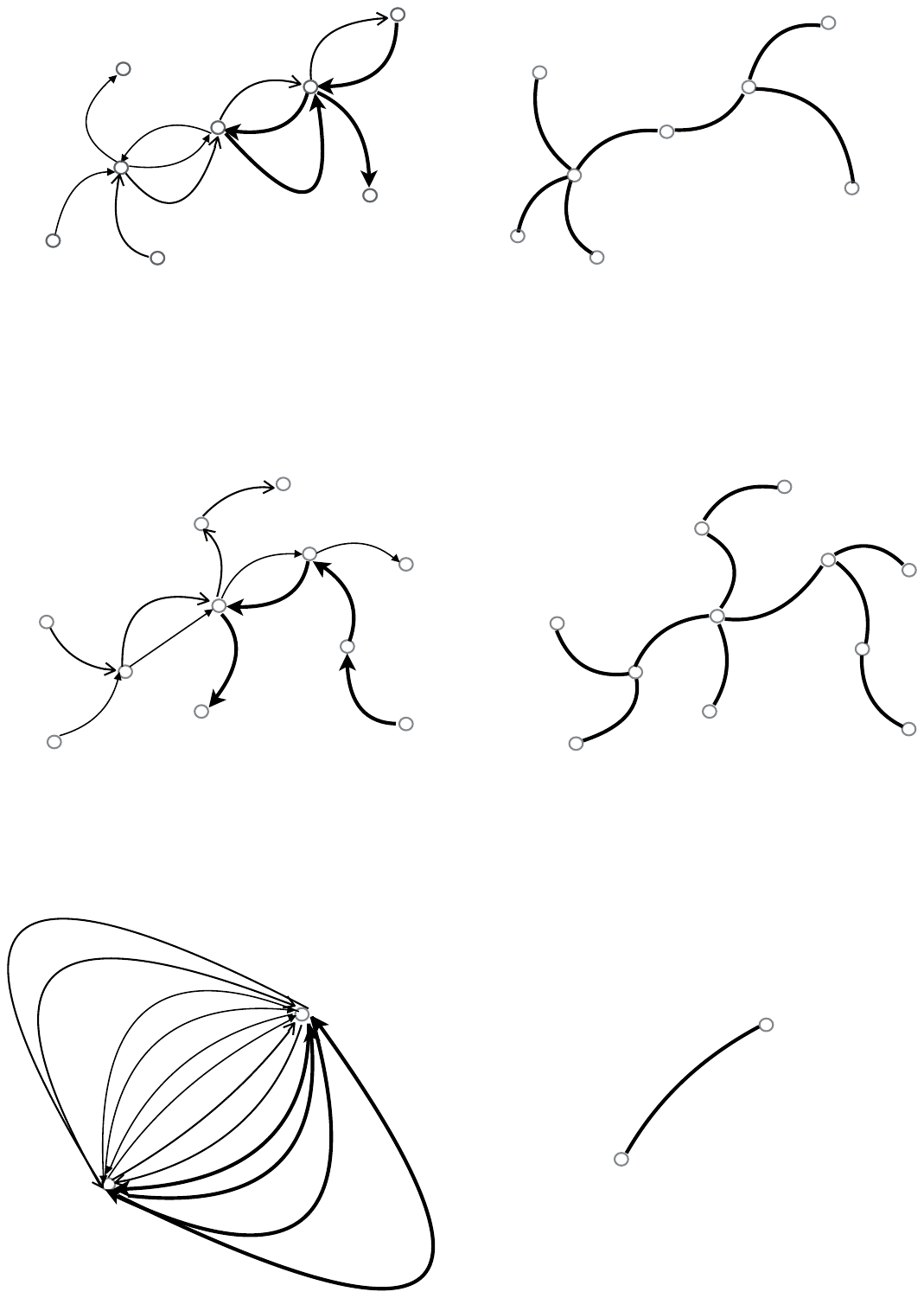}}
\caption{{  {\it Tree-type diagrams $\CT_3^{(4)}$, $\CT_3^{(4, \max)}$, $\CT_3^{(4, \min)}$  and 
their graphs 
}}}
\end{figure}


\vskip 0.2cm 

{\bf Lemma 3.3.} {\it 
For any connected diagram $\hat \CD_k^{(q)}(\pi_s)$ with non-trivial partition $\pi_s$, $s\ge 2$,  where $s$ 
is the number of subsets in $\pi_s$,  the following asymptotic relation holds,
$$
\rW_{N,c,R} \big(\hat \CD_k^{(q)}(\pi_s)\big) = 
o\big( \rW_{N,c,R}(\hat \CD_k^{(q)})(\pi_1)\big), \quad (N,c, R)_0\to\infty.
\eqno (3.18)
$$
}

{\it Proof.} As it is easy to see, any  non-trivial partition produces dotted 
e-arcs in the initial diagram  $\hat \CD_k^{(q)}$  and each of the  dotted e-arcs   adds the factor $c/N$ to the total weight of the diagram. 
Indeed, let us consider an edge $e_{ij}= (\rv_i, \rv_j)$ of $\hat \CD_k^{(q)}$ 
and assume that there are  $f$ e-arcs of $\hat \CD_k^{(q)}$ that became dotted under the action of $\pi_s$. 
We denote this number by $f= f(\pi_s, e_{ij})$. 
Then the weight of the edge $w(e_{ij})$ is given by (cf. (3.10))
$$
w(e_{ij})= \prod_{l=1}^{f(\pi_s, e_{ij})+1}  \ {c\over N}  \, 
h_{\kappa_l(\pi_s, e_{ij})}\left( {x_i-x_j\over R}\right)
\le \left({c\over N}\right)^{f(\pi_s,e_{ij})+1} \ 
 h_{m(v_i,v_j)} \left( {x_i-x_j\over R}\right),
%
$$
where multiplicities $\kappa_l(\pi_s, e_{ij})$ are such that 
$$
 \sum_{l=1}^{f(\pi_s, e_{ij})+1} \kappa_l(\pi_s, e_{ij}) = m(\rv_i,\rv_j). 
$$
Then (3.18) follows.  $\Box$

According to Lemmas 3.1 and  3.3,
we can rewrite relation (3.11) in the following form,
$$
\rcum_k(Y^{(q)}) = \sum_{\hat D_k^{(q)}\in \, 
\fD_k^{(q, \rcon )}} 
\rW_{N,c,R}\big(\hat D_k^{(q)}\big) (1+o(1)),
\quad (N,c,R)_0\to\infty,
\eqno (3.19)
$$
where the sum runs over the set of all connected diagrams 
that we denote by $\fD^{(q, \rcon)}_k$.  
Everywhere below, we  
omit 
 hats in denotations 
$\hat D_k^{(q)}$.


\section{Cumulants  of random variables  $Y$}

Let us consider the first element of the sum (3.5) that corresponds 
to the trivial partition $\pi_1= \{1, 2, \dots, k\}$. In this case there is no dotted e-arcs in $\CD_k^{(q)}(\pi_1)$. 
We say that   type-II-diagram $D_k^{(q)}$ is of tree-type if its graph $G_k^{(q)}= G(D_k^{(q)})$ 
is such that \mbox{$|\CV(G_k^{(q)})|= |\CE(G_k^{(q)})|+1$}. 
We denote such tree-type-II-diagrams by 
$\CT_k^{(q)}(\CL_k^{(q)})$. 
We denote  the set of all tree-type diagrams by $\fT_k^{(q)}$. 
If $\CT_k^{(q)}$ is such that 
$|\CE(G(\CT_k^{(q)}))|= k(q-1)+1$, then we say that 
this tree-type diagram is the maximal one and denote it by 
$\CT^{(q,\max)}_k$. 
We denote the set of all maximal tree-type diagrams by $\fT_k^{(q, \max)}$. 
The number of such diagrams $|\fT_k^{(q,\max)}|=t_k^{(q)}$ is studied in 
\mbox{Section 7.}    
If $\CT_k^{(q)}$ is such that 
$
|\CE(G(\CT_k^{(q)}))|=1,
$
then we say that this  tree-type diagram $\CT^{(q)}_k$ is  the minimal 
one and denote 
 such diagram by $\CT^{(q,\min)}_k$.  We denote by $\fT^{(q,\min)}_k$ the set of all minimal tree-type diagrams.
On Figure 2, we give examples of  tree-type diagram $\CT_3^{(4)}$, 
maximal tree-type diagram $\CT_3^{(4, \max)}$ and minimal tree-type diagram
$\CT_3^{(4, \min)}$ as well as their graphs.

\subsection{Connected diagrams and tree-type diagrams   for $Y$-models}

Given a family of diagrams  $\rD$, we denote its weight by
$$
\CW_{N,c,R} \big(\rD\big)= 
\sum_{ D_k^{(q)} \in \rD}
\rW_{N,c,R}
( D_k^{(q)}),
\eqno (4.1) 
$$
where $W_{N,c,R}(D_k^{(q)})$ is  determined by (3.10).


\vskip 0.2cm 
{\bf Lemma 4.1.} {\it In the limit $(N,c,R)_0\to\infty$ (2.6), 
\vskip 0.1cm 
1) if $cR/N\gg 1$, then 
$$
\CW_{N,c,R} \big(\fD_k^{(q,\rcon)}\big)= 
 \CW_{N,c,R}\big(\fT_{k}^{(q,\max)}\big) (1+o(1)),
\quad (N,c,R)^{(1)}_0\to\infty,
\eqno (4.2)
$$
 where we denoted by $(N,c,R)_0^{(1)}\to\infty$
 the limiting transition (2.6) such that $cR/N\gg 1$;
}

\vskip 0.2cm 
{\it 2) if $cR/N= s$, 
$$
\CW_{N,c,R} \big(\fD_k^{(q,\rcon)}\big)= 
\CW_{N,c,R}(\fT_k^{(q)}) (1+o(1)),
\quad (N,c,R)^{(2)}_0\to\infty,
\eqno (4.3)
$$
where we denoted by $(N,c,R)^{(2)}_0\to\infty$ 
the limiting transition (2.6) such that $cR/N=s$;
}  

\vskip 0.1cm 
{\it 3) if $cR/N\ll 1$, then 
$$
\CW_{N,c,R} \big(\fD_k^{(q,\rcon)}\big)= 
\CW_{N,c,R}(\fT_k^{(q,\min)}) (1+o(1)), \quad 
(N,c,R)^{(3)}_0\to\infty,
\eqno (4.4)
$$
where we denoted by $(N,c,R)^{(3)}_0\to\infty$ 
the limiting transition (2.6) such that $cR/N\ll 1$.
}
\vskip 0.2cm 
{\it Proof.}
In view of   Lemma 3.2, 
 we attribute to each diagram $D_k^{(q)}$ 
its order 
$$ 
\O( D_k^{(q)})= p_N^E N R^{V-1} = \left( {c\over N}\right)^{E} N R^{V-1},  
\eqno (4.5)
$$
where $E=E(D_k^{(q)})= |\CE(G( D_k^{(q)}))|$ and $V=V(D_k^{(q)})= |\CV(G(D_k^{(q)}))|=|\CV(D_k^{(q)})|$. 
Each diagram can be classified according to the value of  
$(E,V)= (E(D_k^{(q)}),V(D_k^{(q)}))$ and placed into corresponding cell (box)
of the plane with the Descartes coordinates. We denote by $\cal S$ the 
collection of all such possible boxes.
$$
\begin{array}{rrrr}
& P \ \ \ &  \ \ \  & \\
\\ 
Q \ \ \ \ & \underline{A}\ \ \   &  \ \ \ \ \ \ V_{\max} \ \ \\
\\ S \ \ \ \ \ \ \  \ \ 
\underline{B}\ \ \  \  &  \underline{K}\ \ \  &  \  V_{\max}-1& \\ \ \ \ \ \ \ \ 
\ \ \ \ 
\underline{L}   \ \ \ \ & \vdots \ \ \ \  &  \vdots \ \ \ \ \, 
\\ 
\underline{F}\ \ \ \ \cdots  \ \  \ \ \ \ \ \ \ \  \vdots \ \ \  \ \,  & \vdots \ \ \ \   &  \ \ \ \ \ \ \ \ \ \ \ 2 \ \ \ \ \\
\\ 
\, 1\ \ \ \ \dots \  \ \ \ E_{\max}-1& \ \ \ E_{\max} \ \ & 
&\\
\end{array}
$$

\begin{figure}[htbp]
\centerline{
}
\caption{{{\it Classification of connected diagrams $ D_k^{(q)}$ on the plane $(E,V)$ }}}
\end{figure}

Let us determine the value of maximally possible number of edges $E_{\max}$ of graphs $G(D_k^{(q)})$. 
In connected diagram  $\CD_k^{(q)}$,
each element $\l'_q$ 
can be attributed by a minimal arc that joins $\l'_q$ with $\l''_q$ that has the minimal number
among those joined to $\l'_q$. 
Then, by recurrence, it is easy to see  
that any  connected diagram $\CD_k^{(q)}$ contains at least $k-1$ arcs.

If $\CD_k^{(q)}$ contains exactly $k-1$ arcs, then 
$|\CV(\CD_k^{(q)})|= k(q+1)- 2(k-1)=k(q-1)+2$
and $|\CE(G(\CD_k^{(q)}))|= kq- (k-1)= k(q-1)+1$. This can be proved again by recurrence. 
By definition,  $G(\CD_k^{(q)})$ is a tree, and this is the maximal tree,
$G(\CD_k^{(q)})= \CT_k^{(q, \max)}$. Thus, the box $A$ of Figure 3 contains 
 maximal tree-type diagrams $\CD_k^{(q)}$  with  $V_{\max}= k(q-1)+2$ and  
 $E_{\max}=k(q-1)+1$. 
 
 All graphs $G=G(\CD_k^{(q)})$ of connected diagrams are, in their turns, connected and therefore 
 the  box $Q= (V_{\max}, E_{\max}-1)$ is empty; the same is true for any box
 with coordinates $(V_{\max}, E_{\max}-p)$ with $p\ge 1$. The same is true for the box $S= (V_{\max}-1, E_{\max}-2)$
 and all other boxes with coordinates $(E',V')$ such that $E'<V'-1$ as well . 
 
We denote  by $D_A$ and $D_K$ diagrams that belong to  boxes $A$ and $K$, respectively.
It follows from (4.5)  that $\O(D_A) = \O(D_K) R$  and therefore 
 $\O(D_K)\ll \O(D_A)$ in the limit (2.6). 
Also we can write that $\O(D_L)\ll \O(D_B)$. 
Let us note that  diagrams of any box $I=(E,V)$ situated under the main diagonal
have the order much smaller than that of diagrams  of 
the corresponding diagonal box $(E,V')$, $V'>V$.

Let us consider the main diagonal of $\CS$
and denote by $J_l$ the boxes with coordinates 
$(k(q-1)+1-l, k(q-1)+2-l)$, $0\le l\le k(q-1)$. 
Since for the graphs $G(D)$ of any diagram $D$ of these boxes $V=E+1$, 
we conclude that $D$ is a tree-type diagram.
It follows from (4.5)  that 
$$
\O(D_A)= \O(D_{J_l}) \times 
\left({c\over N} R\right)^{l} , \quad 
0< l\le k(q-1). 
$$
If $cR\gg N$, then $\O(D_A)\gg \O(D_{J_l})$ for all $l$ and the leading contribution to (4.1) is given by the maximal tree-type diagrams
with the order 
$$
\O(D_A)= \left( {c\over N} \right)^{k(q-1)+1} N R^{k(q-1)+1}
= \left({cR\over N}\right)^{k(q-1)} cR.
\eqno (4.6)
$$
Relation (3.14) means that 
$
\rW_{N,c,R}(D)= \O(D)(1+o(1)), 
\quad (N,c,R)_0^{(1)}\to\infty 
$
and therefore for any diagram $D'\notin A$ we get 
$$
\rW_{N,c,R}(D') = o\left( \left({cR/ N}\right)^{k(q-1)} cR\right),
\quad (N,c,R)_0^{(1)}\to\infty. 
$$
This observation together with relation
$$
\sum_{D\in A} \rW_{N,c,R}(D) = 
\CW_{N,c,R}(\fT_k^{(q,\max)})
$$
implies relation (4.2). On Figure 3, we present an example of the maximal tree-type diagram $\CT_3^{(4)}$ and its graph $G(\CT_3^{(4)})$. 
If $cR/N= s$, then diagrams of the diagonal boxes
$A,B, \dots, F$ 
are all of the order $O(s)$ in the limit $(N,c,R)_0^{(2)}\to\infty$. It is clear that diagrams of 
any diagonal box are the tree-type ones and then (4.3)  follows.
Relation 
$$
\O(D_{J_l})= \O(D_F) \times \left({c\over N} R\right)^{k(q-1)- l},
\quad 0\le l < k(q-1)
$$
shows that n the third asymptotic regime $(N,c,R)^{(3)}\to\infty$,
 the leading contribution to (4.1) is given by minimal
tree-type diagrams from the box $F$ that have one edge of multiplicity $kq$
such that $
\O(D_F)= cR.
$
Then (4.4) follows. Lemma 4.1 is proved.  $\Box$

\vskip 0.2cm

{\bf Lemma 4.2.}  {\it 
In all of the three asymptotic regimes of Lemma 4.1, the following relation is true,
$$
\lim_{(N,c,R)_0\to\infty} {W_{N,c,R}^{(\a)}(\CT_k^{(q)})\over N(p_NR)^{E}  } 
= 
w^{(\a)} (\CT_k^{(q)}),  \quad \a=0,1,
\eqno (4.7)
$$ 
where $E= |\CE(G(\CT_k^{(q)}))|$ and the weight coefficient $w^{(\a)}(\CT_k^{(q)})$ is given by 
$$
w^{(\a)} (\CT_k^{(q)})= 
\int_{-\infty}^{\infty} \cdots \int_{ -\infty} ^{ \infty}
\prod_{
{\stackrel{i,j:}
 {\{v_i,v_j\}\in \CE(\CT_k^{(q)}) }}
} \ h_{m(i,j)}^{(\a)}(x_i-x_j)
\big\vert_{x_1=0} \prod_{l=2}^{V}
\, dx_l
= \prod_{
{\stackrel{i,j:}
 {\{v_i,v_j\}\in \CE(\CT_k^{(q)}) }}
} \ V_{m(i,j)}^{(\a)}
\eqno (4.8)
$$ 
 where $V_m^{(0)}= V_0$ (2.7) and $V_m^{(1)}$ is determined by (2.14). 
 }

\vskip 0.2cm 
{\it Proof.}
Let us consider an auxiliary  tree-type diagram  
$\CT^{(q_1, \dots, q_r)}=\rT_r$ 
constructed with the help of $r$ elements $\l_{q_1}, \dots, \l_{q_r}$
of length $q_i\ge 1$, $1\le i\le r$ such that 
the number of vertices of this diagram is given by 
$V= | \bar q_r|-r+2$, where we denoted 
$\bar q_r=(q_1, \dots, q_r)$ and $|\bar q_r| = q_1+ \dots + q_r$. 
We will also use denotation $\tau_{\bar q}$ for such tree-type diagrams. 
Using  recurrence by $l=|\bar q_r|\ge 1$, we prove the following statement $A_l$,  
$$
\rA_l:\lim_{(N,R)_0\to\infty} {1\over 
N( p_N R)^{E}    } W_{N,R}(\CT^{(q_1, \dots, q_r)})
= w^{(\a)}(\CT^{(q_1, \dots, q_r)}).
\eqno (4.9)
$$

The initial step is given by the diagram $\rT_r=\CT^{(1,1, \dots, 1)}$ that contains one multiple edge $(v_1,v_2)$.  Then according to (3.10), we have 
$$
W_{N,c,R} (\CT^{(1, 1, \dots, 1)})=
{c\over N} \sum_{(s_1, s_2)\in [-n,n]^2} h_{m(1,2)}\left( {s_1-s_2\over R}\right). 
$$ 
To study the limit of this expression, we  perform the following standard  actions:
we restrict the sum over $s_1$ to the sum over interval
$[-n+RL, n-RL]$, then for each given $s_1$ from this interval we replace the 
sum over the interval $[-n,n]$ by the sum 
of $h((s_1-s_2/R))$ over the set $s_2: |s_1-s_2|\le RL$,
such that the result is represented by a  value that does not depend on $s_1$ and is 
close to $\int_{-\infty}^\infty h(t)dt$. To do this, we write 
that 
$$
{1\over cR} W_{N,c, sR} (\CT^{(1, 1, \dots, 1)})=
{1\over N} \sum_{s_1\in [-n+RL, n-RL]} {1\over R} \sum_{s_2\in [-n,n]}
h_{m(1,2)}\left( {s_1-s_2\over R}\right) + \Sigma_1^{(N,R,L)},
\eqno (4.10)
$$
where  
$$
\Sigma_1^{(N,R,L)}= {1\over N} \sum_{s_1\in [-n,n]\setminus I_1} {1\over R} \sum_{s_2\in [-n,n]}
h_{m(1,2)}\left( {s_1-s_2\over R}\right)
$$
and $I_1= [-n+RL, n-RL]$.
Taking into account  elementary upper bound (see  (3.17)),
$$
\sup_{s_1} {1\over R} \sum_{s_2\in [-n,n]} h_{m(1,2)}\left( {s_1-s_2\over R}\right)
\le {1\over R}\big( h_{m(1,2)}( 0) + h_{m(1,2)}(x_1)\big) + 2 \int_0^\infty h(t) dt,
$$
we can conclude   that for any given $L$, the following relation holds, 
$$
\Sigma_1^{(N,R,L)}=  O({2RL/ N})=   o(1), \quad (N,c, R)_0\to\infty.
\eqno (4.11)
$$

Next, 
we  write that
$$
{1\over NR} \sum_{s_1\in I_1}  \sum_{s_2\in [-n,n]}
h_{m(1,2)}\left( {s_1-s_2\over R}\right)
={1\over NR} \sum_{s_1\in I_1}  \sum_{s_2: |s_1-s_2|\le RL}
h_{m(1,2)}\left(  {s_1-s_2\over R}\right)+
\Sigma_2^{(N,R,L)},
$$
where we denoted 
$$
\Sigma_2^{(N,R,L)}={1\over NR} \sum_{s_1\in I_1} 
\sum_{\stackrel {s_2\in [-n,n] }{ |s_1-s_2|> RL}}
h_{m(1,2)}\left(  {s_1-s_2\over R}\right).
$$
Using an elementary upper bound
$$
\sup_{s_1\in [-n,n]} 
{1\over R} 
\sum_{\stackrel {s_2\in [-n,n], }{ |s_1-s_2|> RL}}
h_{m(1,2)}\left(  {s_1-s_2\over R}\right)
\le 
\sup_{s_1\in [-n,n]} 
{1\over R} 
\sum_{\stackrel {{s_2 \in {\mathbb Z} }}{ |s_1-s_2|> RL}}
h_{m(1,2)}\left(  {s_1-s_2\over R}\right)
$$
$$
= {1\over R} 
\sum_{ |s_2| >RL  }
h_{m(1,2)}\left(  {-s_2\over R}\right) \le \int_{-\infty}^{-L} h_{m(1,2)}(t) dt
+ \int_{L}^\infty h_{m(1,2)}(t) dt= \vep_1^{(L)},
$$
we conclude that 
$$
\Sigma_2^{(N,R,L)}=  O( \vep_1^{(L)}/N),
\eqno (4.12)
$$
where 
$\vep_1^{(L)}$ tends to zero as $L\to\infty$. 
Using two elementary relations,
$$
{1\over R} \sum_{s_2: |s_1-s_2|\le RL}
h_{m(1,2)}\left(  {s_1-s_2\over R}\right) - \int_{-L}^L 
h_{m(1,2)}(t) dt = O(1/R)
\eqno (4.13)
$$
and 
$$
\int_{-\infty}^\infty
h_{m(1,2)}(t) dt - 
\int_{-L}^L 
h_{m(1,2)}(t) dt = \vep_2^{(L)},
\eqno (4.14)
$$
where $\vep_2^{(L)} $ tends to zero as $L\to\infty$,
we deduce from (4.10), (4.11) and (4.12)  that 
$$
{1\over cR} W_{N,c, R} (\CT^{(1,1, \dots, 1)})
= 
\int_{-\infty}^\infty
h_{m(1,2)}(t) dt  + O( \vep_1^{(L)}/N+ \vep_2^{(L)}
+ 1/R)+ o(1). 
$$
Since $L$ can be chosen arbitrary large, this relation 
shows  that the statement $\rA_1$ of (4.9) is true. 

\vskip 0.1 cm

Let us consider the general case $\rA_l$. 
In  tree-type diagram $\CT^{(q_1, \dots, q_r)}= \rT_r$, we consider the set of  extreme vertices of leafs 
and determine  the maximal vertex of this set; we denote its number 
by $\kappa$. We denote by $\chi$ the number of the vertex 
$v_\chi$ such that 
$\{v_\chi, v_{\k }\} \in \CE(\rT_r)$. 
Finally, we denote by 
$
[-n,n]^l_{(\chi, R,L)}
$
the set such that the interval number $\chi$ is restricted to 
the interval  $[-n + RL, n-RL]$. Then we can write that
$$
{W_{N,c,R} ( \rT_r)\over N (p_N R)^E} =
\sum_{
(s_1, \dots , s_l) \in [-n,n]^l_{(\chi, R,L)}}\ 
\prod_{
{\stackrel{1\le i<j\le l:}
 {\{v_i,v_j\}\in \CE(\rT_r)}}
}
 \ h_{m(i,j)}\left({s_i-s_j\over R}\right)
+ \tilde \Sigma_1^{(N,R,L)}(\chi),
\eqno (4.15)
$$
where  $l$ is the total number of vertices of $\rT_r$ and 
$$
\tilde \Sigma_1^{(N,R,L)}(\chi)=
 {1\over R^E N}
\sum_{(s_1, \dots , s_l) \in [-n,n]^l\setminus [-n,n]^l_{(\chi, R,L)}}\ 
\prod_{
{\stackrel{1\le i<j\le l:}
 {\{v_i,v_j\}\in \CE(\rT_r)}}
} \ h_{m(i,j)}\left({s_i-s_j\over R}\right).
$$
In this sum, we consider  the vertex $v_\chi$ as the root one and  
attribute the normalizing factor $1/N$ to the sum over $s_\chi\in [-n + RL, n-RL]$.
It remains to show that for any given value of $s_\chi$, 
the sum over variables $\{s_1, \dots, s_l\} \setminus \{s_\chi\}$
multiplied by $1/R$ is bounded from above. 
This can be done by recurrence with the help of (4.13)
and we omit here the elementary reasoning.  Then we can write,  as before, an asymptotic relation
 $
 \tilde  \Sigma_1^{(N,R,L)}(\chi)=   O(2RL/N)=  o(1)$,  
 $(N,c,R)_0\to\infty.
$ 
Regarding the first term of the right-hand side of (4.15), we observe that it  contains the factor
$$
{1\over N} \sum_{s_\chi \in [-n + RL, n-RL]} \ {1\over R} \sum_{s_{\k }\in [-n,n]}
 \ h_{m(\chi,{\k })}\left({s_\chi-s_{\k }\over R}\right)
 $$
that we treat exactly as it is done in the proof of (4.12). Then, using 
relations (4.13) and (4.14), we obtain that 
$$
{W_{N,c,R} ( \rT_r)\over N (p_NR)^E}  = 
\int_{-\infty}^\infty
h_{m(\chi,{\k })}(t)dt (1+o(1))
$$
$$
\times
\  \sum_{(s_1, \dots, s_l)\setminus 
\{s_{\k}\} \in [-n,n]^{l-1}_{\chi, R,L}}
\ \prod_{
{\stackrel{1\le i<j\le l-1:}
 {\{v_i,v_j\}\in \CE(\rT_{r'})}}
}
 \ h_{m(i,j)}\left({s_i-s_j\over R}\right),
 \quad (N,c,R)_0\to\infty,
 \eqno (4.16)
$$
where $\rT_{r'}$ is obtained from $\rT_r$ 
by erasing the multi-edge $\{v_\chi, v_{\k }\}$. 
In the last factor  of (4.16), it remains to pass from the sum indicated
to the summation of values $(s_1, \dots, s_l)\setminus 
\{s_{\k }\} $ over $[-n,n]^{l-1}$ such that $\rA_{l-1}$ of (4.9) can be used. This transition can be justified by the same reasoning as used 
in the  study of  $ \Sigma_1^{(N,R,L)}(\chi)$. We omit these elementary arguments. 
Relation (4.9) is proved. The first equality of (4.8) is an 
obvious consequence of (4.9). The proof of the second equality of (4.8) is also elementary.
Lemma 4.2 is proved. $\Box$

\subsection{Proof of Theorems 2.1 and 2.2}

We start with  the first limiting transition  
$(N,c,R)_0^{(1)}\to\infty$ when $cR\gg N$. 
It follows from  (3.19) and  relation (4.2) of  Lemma 4.1 
that 
$$
\rcum_k(Y^{(\a,q)})= \CW_{N,c,R}^{(\a)}(\fT_k^{(q,\max)})
(1+o(1))
,
\quad (N,c,R)_0^{(1)}\to\infty.
\eqno (4.17)
$$
Relation (4.6) shows that all maximal tree-type diagrams
are of the same order of magnitude
$cR (cR/N)^{k(q-1)}$. Then  it follows from relation (4.7) 
of Lemma 4.2  that 
$$
\rcum_k(Y^{(\a,q)}) = 
cR \left({cR\over N}\right)^{k(q-1)}
\sum_{\CT_k^{(q)}\in \, \fT_k^{(q,\max)}} 
w^{(\a)}(\CT_k^{(q)}) (1+o(1)),
\quad (N,c,R)_0^{(1)}\to\infty
$$
and finally, that
$$
\lim_{(N,c,R)_0^{(1)}\to\infty} {1\over cR} 
\rcum_k
\left(  \hat Y^{(\a,q)}\right)= 
\begin{cases}
 \Phi_k^{(q,1)},  &  \text{if  $\a=0$}  , \\
\Xi_k^{(q,1)},  \quad   & 
\text {if $\a=1$,
}
\end{cases}
\eqno (4.18)
$$
where $\hat Y^{(\a,q)}= (N/cR)^{q-1} Y^{(\a,q)}$ and where 
$$
\Phi_k^{(1)}=\sum_{\CT_k^{(q)}\in \fT_k^{(q,\max)}} 
w^{(0)}(\CT_k^{(q)})
= | \fT_k^{(q,\max)}| =t_k^{(q)}
\eqno (4.19)
$$
represents the total number of the maximal tree-type diagrams 
constructed with the help of $\l_q$-elements and 
$$
 \Xi_k^{(q,1)}=
\sum_{\CT_k^{(q)}\in \fT_k^{(q,\max)}} 
w^{(1)}(\CT_k^{(q)}),
\eqno (4.20)
$$
is the total sum of weighted maximal tree-type diagrams, 
with   formula (4.8) 
used for the weight  coefficient $w^{(a)}(\CT_k^{(q)})$.
Relation (4.17) proves existence of limits (2.8) and (2.15). 
Explicit expressions for the limiting values
$\Phi_k^{(q,1)}$ and $\Xi_k^{(q,1)}$ given by relations (2.11)  and (2.18)
 will be obtained  in Section 7.

Let us consider the second asymptotic regime when 
$cR/N= s$ in the limiting transition $(N,c,R)_0\to\infty$ (2.6).
In this case, relations (4.17) and  (4.18)
take form
$$
\rcum_k(Y^{(\a,q)})= 
\CW_{N,c,R}(\fT_k^{(q)})
(1+o(1)),
\quad (N,c,R)_0^{(2)}\to\infty,
\eqno (4.21)
$$
$$
\lim_{(N,c,R)_0^{(2)}\to\infty} {1\over cR} 
\rcum_k\left(
Y^{(\a,q)}\right)= 
\begin{cases}
 \Phi_k^{(q,2)}(s),  &  \text{if  $\a=0$}  , \\
\Xi_k^{(q,2)}(s),  \quad   & 
\text {if $\a=1$,
}
\end{cases}
\eqno (4.22)
$$
where 
$$
\Phi_k^{(2)}(s)= \sum_{s=1}^{k(q-1)+1} 
\sum_{
\stackrel{q_1\ge 1, \dots q_k\ge 1} {q_1+ \dots +q_k=s} } 
| \fT^{(q_1,\dots, q_k)}|
$$
and $\fT^{(q_1,\dots, q_k)}$ is the set of all tree-type diagrams
constructed with the help of elements $\l_{q_1}, \dots, \l_{q_k}$. 
Cardinality of these families of trees will be considered
and relations (2.12) and (2.13) will be obtained in Section 7. 
Explicit expression for $\Xi_k^{(q,2)}$ of (4.22) is rather complicated and we do not present it here.

In the third asymptotic regime, 
we have convergence
$$
\lim_{(N,c,R)_0^{(3)}\to\infty} {1\over cR} 
\rcum_k\left(
Y^{(\a,q)}\right)= 
\begin{cases}
 \Phi_k^{(q,3)}(s),  &  \text{if  $\a=0$}  , \\
\Xi_k^{(q,3)}(s),  \quad   & 
\text {if $\a=1$,
}
\end{cases}
\eqno (4.23)
$$
where, according to (2.14) and (3.10), the term 
$\Phi_k^{(q,3)}$ is given by $V_0$ multiplied 
by the number {$| \fT^{(q,\min)}|=t_k^{(\min)}$}
and $\Xi_k^{(q,3)}$ is given by $V_{kq}$ multiplied by $t_k^{(\min)}$.
It is easy to see that $t_k^{(\min)}$ represents the number of ways to put
$k-1$ oriented edges on the first oriented edge that is equal to $2^{k-1}$. 
Taking into account  this observation, we deduce from (4.23)  relations (2.10) and (2.17). Theorems 2.1 and 2.2 are proved. $\Box$

\section{Cumulants of random variables $X^{(3)}$}

\begin{figure}[htbp]
\centerline{\includegraphics
[height=15cm, width=12cm]{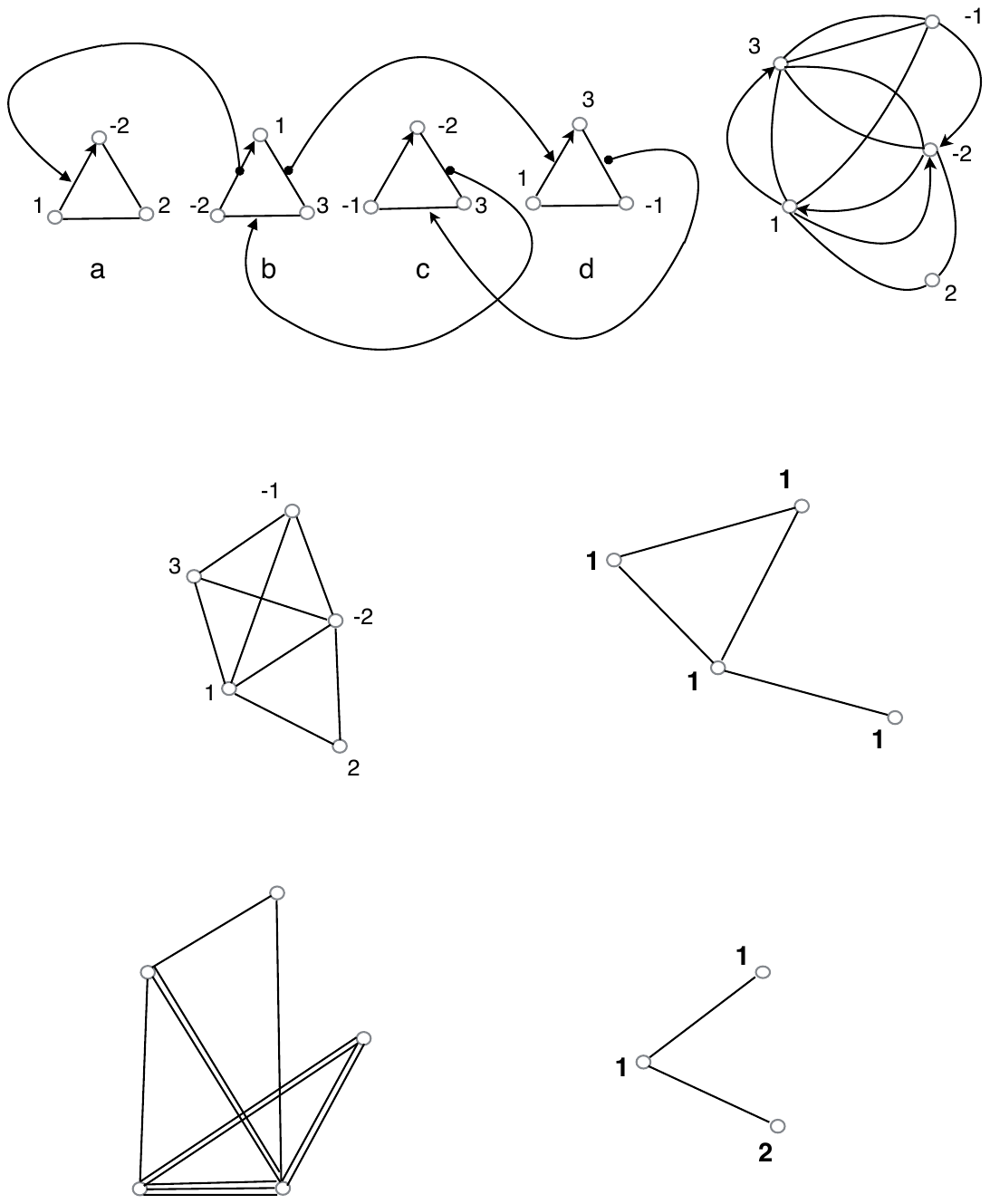}}
\caption{{  {\it Diagrams $\CD_4^{(3)}$, $D_4^{(3)}$, graphs
$G(D_4^{(3)})$,  $\CG(D_4^{(3)})$; diagram $\tau_4^{(\mu)}$ and $\CG(\tau_4^{(\mu)})$ 
}}}
\end{figure}

In this section, we consider asymptotic behavior of cumulants of random variables 
$$
 X^{(\a, 3)}_{N,c,R}= \sum_{i_1, i_2, i_3 \in \rL_N } \ 
\left( A^{(\a)}_N\right)_{i_1 i_2} \left( A^{(\a)}_N\right)_{i_2 i_3}  
\left(  A^{(\a)}_N\right)_{i_3 i_{1}}
=\sum_{i_1, i_2, i_3 \in \rL_N } \ \CX^{(\a,3)}(\la i\ra_3),
\eqno (5.1)
$$
where we denoted $\la i\ra_3= (i_1, i_2,  i_3)$. 
As in Section 3, we can write  the semi-invariant representation
$$
\rcum_k(  X^{(\a, 3)}_{N,c,R})= 
\sum_{\CL_k^{(q)}(N) \in \rL^{kq}_N}
 cum_k\left\{ \la   \CX^{(\a, 3)}\ra_1, \dots, \la   \CX^{(\a,3)}\ra_k\right\},
 \eqno (5.2)
$$
where $\la \CX^{(a,3)}\ra_j= \CX^{(\a,3)}(\la i\ra^{(j)}_3)$,
$
\la i\ra^{(j)}_3= (i_1^{(j)}, i^{(j)}_2,i^{(j)}_3)$, 
and 
$
\CL_k^{(3)}(N)= \left( \la i\ra^{(1)}_{3}, \   \dots \  , 
\la i\ra^{(k)}_{3}\right)_N
$ is a realization of the set of variables $\{ \la i\ra^{(1)}_{3}, \   \dots \  , 
\la i\ra^{(k)}_{3}\}$ (cf. (3.8)).

To study  mixed cumulants 
$cum_k\left\{ \la   \CX^{(\a, 3)}\ra_1, \dots, \la   \CX^{(\a,3)}\ra_k\right\}$, 
we  repeat  considerations of sub-section 3.2 
 with the only  difference that   $\l_q$-elements are replaced
 by $\mu_3$-elements.  
In this case, diagrams  are constructed with the help of 
triangle elements $\mu_3$ with cyclically oriented edges. For each $\mu$-element, we indicate
 orientation of one edge only. We denote these diagrams by $\CD_k^{(\mu)}$ to distinguish them from 
 $\CD_k^{(q)}$ used in Section 3. 

Rigorous definition of $\CD_k^{(\mu)}$ is similar to that of $\CD_k^{(q)}$: given a realization $\CL_k^{(3)}(N)$,
we construct a type-I-diagram $\CD_k^{(\mu)}= \CD(\CL_k^{(3)}(N))$ as the set of $k$ $\mu$-elements, where  
positively or negatively coupled edges of $\mu$-elements are joined by positively or negatively oriented 
e-arcs;  in $\CD_k^{(\mu)}$ the vertices  attributed by the same values of $\rL_N$
and  not attached to the edges joined by e-arcs are joined by non-oriented v-arcs. We minimize 
the collection of e-arcs by keeping the e-arcs that join nearest $\mu$-elements and erasing
all other e-arcs; the set of v-arcs is minimized by the same procedure. 

The type-II-diagrams $D_k^{(\mu)}= D(\CL_k^{(3)}(N))$ are obtained from $\CD_k^{(\mu)}= \CD(\CL_k^{(3)}(N))$ by identifying  vertices   attributed by the same values of $\rL_N$. 
One can  consider $\CD^{(\mu)}_k$ as a multi-graph with ordered vertices. Regarding an edge 
$e=(\rv',\rv'')\in \CE(\CD^{(\mu)}_k)$, we denote by $m(e)$ its multiplicity that counts the simple oriented edges joining 
$\rv'$ and $\rv''$ in both directions.

We determine the weight $W^{(\a,3)}_{N,c,R}(D_k^{(\mu)})$ by equality (3.9) with   subscripts $q$ replaced 
\mbox{by $\mu$}. 
Formulas (3.10), (3.11), (3.12), definitions of connected diagrams $\hat D_k^{(\mu)}$ and non-connected diagrams 
$\ddot D_k^{(\mu)}$, as well as the statements and the proofs of Lemma 3.1, Lemma 3.2 and Lemma 3.3  remain  valid with the same obvious replacement of subscripts $q$ \mbox{by $\mu$.}  
Similarly to (4.1), we introduce the sum of  weights of connected diagrams,
$$
\CW_{N,R} \big(\fD_k^{(\mu,\rcon)}\big)= 
\sum_{ \hat D_k^{(\mu)} \in \fD_k^{(\mu,\rcon)}}
\rW_{N,R}
( \hat D_k^{(\mu)}),
\eqno (5.3) 
$$
where $\fD_k^{(\mu,\rcon)}$ is the set of all connected 
diagrams constructed with the help of $k$ \mbox{$\mu_3$-elements.}
In what follows, we consider connected diagrams only and omit hats in their denotations.


\subsection{Tree-type diagrams for $X^{(3)}$-model}

We start with the definition of  tree-type diagrams for $X^{(3)}$-model. 
To do this,  we introduce 
 an auxiliary graph $\CG_k= \CG(\CD_k^{(\mu)})$ as follows. Let us  consider $k$ ordered 
 elements 
 $\mu_j, 1\le j\le k$ of the diagram 
 $\CD_k^{(\mu)}$.  If an element $\mu'$ is such that there is no other element $\mu''$
 joined with $\mu'$ by three arcs, we represent it by  a vertex $v'$ attributed by  $1$
 and say that $v'$ is of multiplicity one. 
 If there is a subset of $l$ elements $M_l=\{\mu_{i_1}, \dots , \mu_{i_l}\}$ 
 such that each element of $M_l$ is joined with another element of $M_l$ by three e-arcs,
 we represent this subset by a vertex $v''$ attributed by $l$ and say that $v''$ is of multiplicity $l$.
 Repeating this procedure by recurrence, we get an ordered set of $r$ vertices $v_1, \dots, v_p$
 such that the sum of their multiplicities is equal to $k$. This set determines the set of vertices $\CV(\CG)$. 
 If two vertices $v_i,v_j$ of $\CV(\CG)$ are such that corresponding $\mu$-elements are joined in $\CD_k^{(\mu)}$  
 by an e-arc, 
 we draw an edge $e(v_i,v_j)\in \CE(\CG)$. The graph $\CG_k^{(q)}= \CG(\CD_k^{(\mu)})$ can be considered as 
 the dual one with respect to the multigraph (diagram) $\CD_k^{(\mu)}$. 
Clearly,  $\CD_k^{(\mu)}$ is connected if and only if $\CG(\CD_k^{(\mu)})$ is connected. 
 On \mbox{Figure 4,} we present examples of diagrams $\CD_4^{(3)}$, $D_4^{(3)}$, their graph 
$G(D_4^{(3)})$ and the dual graph $\CG(D_4^{(3)})$ with vertex multiplicities in boldface.

We say that $D_k^{(\mu)}$ is a tree-type diagram
if the dual graph $\CG(D_k^{(\mu)})$ is a tree. 
 We denote  tree-type diagrams by $\tau^{(\mu)}_k$ or simply by $\tau_k$. We say that 
 a  tree-type diagram $\tau_k$ is the maximal one when 
$|\CE(\tau_k)| = 2k+1$ and 
$|\CV(\tau_k)|= k+2$. We denote the maximal tree-type diagrams by 
$\tau_k^{( \max)}$. We say that  a tree-type diagram is the minimal one if
its dual graph 
$\CG(\tau_k^{(\mu,\min)})$ consists of one vertex of multiplicity $k$.
We denote by $\fT^{(\mu,\max)}_k$,  
$\fT^{(\mu)}_k$ and $\fT^{(\mu,\min)}_k$
the ensembles of the \mbox{maximal},  tree-type \mbox{$\mu$-diagrams},
the minimal tree-type $\mu$-diagrams.

{\bf Lemma 5.1.} {\it In the limit $(N,c,R)_0\to\infty$ (2.6),  
\vskip 0.1cm 
i)  if $c^2R/N^2\gg 1$, then 
$$
\CW_{N,c,R} \big(\fD_k^{(\mu,\rcon)}\big) 
 = 
 \CW_{N,c,R} \big(\fT_k^{(\mu, \max)}\big)(1+o(1)),
 \quad (N,c,R)^{(\ri)}_0\to\infty, 
\eqno (5.3)
$$
 where we denoted by $(N,c,R)^{(\ri)}_0\to\infty$ 
the limiting transition (2.6) such that \mbox{$c^2R/N^2\gg 1$;
}}

\vskip 0.2cm 
{\it ii) if $c^2R/N^2= s$,  then 
$$
\CW_{N,c,R} \big(\fD_k^{(\mu,\rcon)}\big) 
 = 
 \CW_{N,c,R} \big(\fT_k^{(\mu)}\big)(1+o(1)),
\quad (N,c,R)^{(\rii)}_0\to\infty,
\eqno (5.4)
$$
where we denoted by $(N,c,R)^{(\rii)}_0\to\infty$ 
the limiting transition (2.6) such that $c^2R/ N^2= s$;
}  

\vskip 0.1cm 
{\it iii)  if $c^2R/N^2 \ll 1$, then 
$$
\CW_{N,c,R} \big(\fD_k^{(\mu,\rcon)}\big) 
 = 
 \CW_{N,c,R} \big(\fT_k^{(\mu, \min)}\big)(1+o(1)),
\quad (N,c,R)^{(\riii)}_0\to\infty,
\eqno (5.5)
$$
where we denoted by $(N,c,R)^{(\riii)}_0\to\infty$ 
the limiting transition (2.6) such that $c^2R/N^2\ll  1$.
}

{\it Proof.  } 
Similarly to Lemma 3.2, we can prove that the weight of any  diagram 
$D_k^{(\mu)}$, in the limit $(N,c, R)_0\to\infty$ (2.6), is given by the following asymptotic relation, 
$$
\rW_{N,c,R} ( D_k^{(\mu)})=
 O\left({c^{E } R^{V -1}\over N^{E-1}}\right), \ \  {\hbox{where}} \ 
 E=E(D_k^{(\mu)})= | \CE( G(D_k^{(\mu)}))| \ {\hbox{and}} \   
V =
|\CV( D_k^{(\mu)})|.
\eqno  (5.6)
$$
We attribute  to  $D_k^{(\mu)}$  expression 
$\O(D^{(\mu)}_k)= NR^{V-1}(c/N)^{E}$  (cf. (4.5)). 
and  distribute  the set  of all possible  connected diagrams 
into cells (boxes) with labels $(E(D_k^{(\mu)}),V(D_k^{(\mu)}))$, with $3\le E\le E_{\max}$ on the plane with
 "coordinates" $(E,V)$. 
If $\CD_k^{(\mu)}$ contains exactly $k-1$ arcs, then $E(D_k^{(\mu)})= E_{\max}= 3k - (k-1)= 2k+1$ and 
$V(D_k^{(\mu)}= V_{\max}= k+2$. This can be proved by recurrence with respect to the maximal leave.
Then clearly the dual graph $\CG(D_k^{(\mu)})$ is a tree and this is a maximal tree.
We denote by $A$ the box with coordinates $(E_{\max}, V_{\max})$. 
Inversely, it can be proved by recurrence that  any diagram $D_k^{(\mu)}\in A$ is such that its dual graph $\CG(D_k^{(\mu)})$ is a maximal tree. 




$$
\begin{array}{rrrr}
& Q \ \ \ & P \ \ \  & \\
S\ \ \ & R \ \ \ & \underline{A}\ \ \   & \ V_{\max}\ \ \\
\underline{B}\ \ \ &\underline{L}\ \ \  & \underline{K}\ \ \   & \ \ V_{\max}-1\\
\vdots \ \ \  & \vdots \ \ \  & \vdots \ \ \  &  \\
\ \ E_{\max}-2& \ \ E_{\max}-1& \ \ E_{\max} &\\
\end{array}
$$
\begin{figure}[htbp]
\centerline{
}
\caption{{ {\it Classification of connected diagrams $ D_k^{(\mu)}$ on the plane 
$(E,V)$}
 }}
\end{figure}
\vskip 0.1cm
\noindent 
\newpage On Figure 5, we show position of the maximal box A with respect to other boxes. 
We 
say that boxes with coordinates $(E_{\max}-2r, V_{\max}-r)$ 
represent the main "diagonal".


By definition, boxes $P$ and $Q$ shown on  Figure  5 are empty, as well as all other boxes with coordinates 
$(E, V_{\max}+k)$, $k\ge 1$. 
Let us show that  $R$ is empty. Any diagram $\CD_k^{(\mu)}\in R$  has at least $k$ e-arcs. 
Then there exists at least one e-arc that can be removed from $\CD_k^{(\mu)}$ without breaking the connectivity. Since 
each e-arc identifies at least two vertices of $\CD_k^{(\mu)}$, then the  diagram obtained from $\CD_k^{(\mu)}$ by removal of the extra e-arc has to fill the boxes of the line $P,Q,...$ that is impossible. 
Therefore,  $R$ contains no diagrams. The same is true for  $S$, as well as for all  boxes with coordinates 
$(E_{\max}-j, V_{\max})$, $j>0$. 


Let us look at the box $B$. If $\CD_k^{(\mu)}$ contains $l$ e-arcs,  then  $E(D_k^{(\mu)})= 3k-l $. Since
$D_k^{(\mu)}\in R$, then $3k-l= 2k-1$ and $l=k+1$. Each e-arc identifies either one or two vertices of $\CD_k^{(\mu)}$. Denoting 
by $l'$ and $l''$ the number of such e-arcs, we obtain relation $3k-l'-2l''= k-1$. Then $l=l'+l''$ gives 
$l''=k$ and $l'=1$. This last e-arc $\fe'$  that identifies edges $(v_1, v_2)$ and $(v_3, v_4)$ identifies only one vertex.
This means that these edges have already a couple of vertices already  identified by some of $l''$ e-arcs. 
Then $\fe'$ is such that it join edges of triangles that give a multiple vertex in the dual graph $\CG(D_k^{(\mu)})$. 
Now it is easy to prove that any diagram $\D_k^{(\mu)}$ of $B$ is a tree-type diagram and that 
 $\CG(D_k^{(\mu)}) $ is a tree with $k-2$ simple vertices and one vertex of multiplicity $2$.
 At the bottom of Figure 4, we depict a tree-type diagram $D_4^{(\mu)}=\tau^{(\mu)}_4$ and its dual graph 
 $\CG(\tau_4^{(\mu)})$ with  multiplicities of vertices  in bold.

Now we proceed by recurrence till the last box  of the main  diagonal (not shown on \mbox{Figure 5}) 
with  coordinates $(3,3)$.  It contains  minimal tree-type diagrams 
represented by one triangle where each edge has multiplicity $k$. 
We say that the main diagonal together with the boxes 
with coordinates $(E_{max}-1 - 2l, V_{\max}-1-l)$ represents the upper border of the set of boxes 
$\CS= \CS_k^{(\mu)}$. 
On Figure 5, we underline the boxes that are not empty. Also it is clear that the boxes below the upper border contain diagrams of much smaller order than diagrams of boxes of the upper border of $\CS$.

Let us consider boxes of the  main diagonal  $I_l= (2(k-l)+1, k+2-l)$, $0\le l\le k-1$ and boxes of the upper border $M_j=(2(k-j), k+1-j)$, $0\le j\le k-2$. It 
follows from (4.14) that 
$$
\Omega (D_{I_l})= \Omega (D_A) \times 
\left( {N^2\over c^2R}\right)^{l},
\quad 1\le l\le k-1
\eqno (5.7)
$$
and 
$$
\Omega (D_{M_j})= \Omega (D_A) \times {N\over cR} \left( {N^2\over c^2R}\right)^{j},
\quad 1\le j\le k-2.
$$
  

 In the first asymptotic regime $(N,c,R)_0^{(\ri)}\to\infty$, relation $c^2 R/N^2\gg 1$ obviously  implies relation $cR/N\gg 1$. Then
 $\Omega(D_{I_l}) \ll \O (D_A)$ for all $1\le l\le k-1$  and $\O (D_{M_j}) \ll \O(D_A)$ for all $1\le j\le k-2$. We see that  the leading contribution to the right-hand side of (5.2) is given by diagrams  
 $D\in A$.  It is clear that these are the maximal 
 tree-type diagrams constructed with the help of $k$ elements $\mu_3$
 and 
 $$
 \O(D_A)= {c^{2k+1} R^{k+1}\over N^{2k}}\ .
 \eqno (5.8)
 $$ 
 Relation (3.14) means that  that 
$
\rW_{N,R}(D)= \O(D)(1+o(1)), 
(N,c,R)_0^{(\ri)}\to\infty 
$
and therefore for any diagram $D'\notin A$ we have
an asymptotic estimate
$$
\rW_{N,R}(D') = o\left( 
\left(c^{2}R/N^2\right)^{k} cR\right),
\quad (N,c,R)_0^{(\ri)}\to\infty. 
$$
This observation 
implies relation (5.3).

 In the second asymptotic regime $(N,c,R)_0^{(\rii)}\to\infty$, relation 
  $c^2 R/N^2= s$ implies  \mbox{$cR/N\gg 1$} again. 
  The diagrams of diagonal boxes $I_l$ are of the same order of magnitude
  $$
  \O(D_{I_l})= cR s^{k-l},
  \eqno (5.9)
  $$
  while  \mbox{$\Omega(D_{M_j})= o(\Omega(D_A))$.} Remembering that diagrams  $D\in I_l$  
 are the tree-type ones, we get asymptotic equality (5.4). 
 On Figure 5, we present an example of diagram $\tau_4^{(\mu)}$ as well as its dual graph 
 $\CG(\tau_4^{(\mu)})$.

In the third asymptotic regime $(N,c,R)_0^{(\riii)}\to \infty$,
$$
\O( D_{I_l})= \O(D_F) \times \left({c^2R\over N^2}\right)^{l}, \quad 1\le l\le k-2
$$
 and 
 $$
 \O (D_{M_j}) = \O(D_F) \times {c\over N} 
 \left( {c^2 R\over N^2}
 \right)^j, \quad 0\le j\le k-2.
 $$
The leading contribution to the right-hand side of (5.2)
 is given by  minimal diagrams $\tau^{(\mu, \min)}_k$ that belong to $F$
 such that 
 $
 \O(D_F)= {c^3 R^2/  N^2}.
 $
Then (5.5) follows. Lemma 5.1 is proved. $\Box$

{\bf Lemma 5.2.} 
{\it  If $\tau_k^{(\mu)}$ is a tree-type diagram, then in all of the three asymptotic regimes of Lemma 5.1  
$$
\lim_{(N,c,R)_0\to\infty} {W_{N,R}^{(\a)}(\tau_k^{(\mu)})\over   N R^{V-1} p_N^E} 
= 
\int_{-\infty}^{\infty} \cdots \int_{ -\infty} ^{ \infty}
\ \prod_{
{{i,j:}
 {\{v_i,v_j\}\in \CE\big(\tau_k^{(\mu)}\big)}}
} \ h_{m(i,j)}^{(\a)}(x_i-x_j)
\big\vert_{x_1=0} \ \prod_{l=2}^V
\, dx_l,
\eqno (5.10)
$$ 
where 
$E$ and $V$ are  the numbers of  edges and  vertices of  the graph $G(\tau_k^{(q)})$, respectively
and $h^{(\a)}_m(x)$ is the function determined by (2.24). }

\vskip 0.2cm 
One can prove Lemma 5.2 by 
using relation similar to (3.10) and by the same recurrence with respect to  the maximal leafs of the graph $\CG_k$, as it is done in the proof of Lemma 4.2. We omit the detailed proof of this  recurrence.

\vskip 0.2cm

\subsection{Proof of Theorem 2.3}

We start with the first limiting transition $(N,c,R)^{(\ri)}_0\to\infty$. We rewrite relation (3.19)  as 
$$
\rcum_k(X^{(\a,3)})= 
\CW_{N,R}^{(\a)}(\fD_l^{(\mu, \rcon)} )
(1+o(1)), \quad (N,c,R)_0^{(\ri)}\to\infty
$$
and deduce with the help of  (5.8) the following asymptotic 
relation,  
$$\rcum_k(X^{(\a,3)})= 
cR \left( {c^2R\over N^2}\right)^{k} \sum_{\tau_k^{(\mu)}\in \fT^{(\mu, \max)}_k}
w^{(\a)}(\tau_k^{(\mu)})(1+o(1)),
\quad (N,c,R)_0^{(\ri)}\to\infty.
$$
Then we can write that 
$$
\lim_{(N,c,R)_0^{(\ri)}\to\infty } {1\over cR} {\rm Cum}_k\left( { N^2\over c^2R} X^{(\a,3)}\right) = 
\Theta_k^{(\a, \ri)}
=\sum_{\tau ^{(\mu)}_k\in \fT_k^{(\mu,\max )}}
w^{(\a)}\big(\tau_k^{(\mu)}\big),
$$
where according to (5.10),
$$
w^{(\a)}\big(\tau_k^{(\mu)}\big) = 
\int_{-\infty}^{\infty} \cdots 
\int_{ -\infty} ^{ \infty}
\ \prod_{
{{i,j:}
 {\{v_i,v_j\}\in \CE\big(\tau_k^{(\mu)}\big)}}
} \ h_{m(i,j)}^{(\a)}(x_i-x_j)
\big\vert_{x_1=0} \ \prod_{l=2}^{V}
\, dx_l,
\eqno (5.11)
$$ 
and $V$ denotes the number of vertices in the tree-type diagram
$\tau_k^{(q)}$. 
This proves  (2.21). 

In the second asymptotic regime $(N,c,R)_0^{(\rii)}\to\infty$,
relation (5.4) impies that 
$$
\rcum_k(X^{(\a,3)})= 
\CW_{N,R}(\fT_k^{(\mu)})
(1+o(1)),
\quad (N,c,R)_0^{(\rii)}\to\infty.
\eqno (5.12)
$$
Then 
$$
\lim_{(N,c,R)_0^{(\rii)}\to\infty } 
{1\over cR} \rcum_k(W^{(\a,3)})
= \Theta_k^{(\a, \rii)}(s),
$$
and convergence (2.22) follows,
where according to (5.9), 
$$
\Theta_k^{(\a, \rii)}(s)
=s^{k-l} \sum_{\tau ^{(\mu)}_k\in \fT_k^{(\mu)}(l)}
w^{(\a)}\big(\tau_k^{(\mu)}\big),
\eqno (5.13)
$$
 where $\fT_k^{(\mu)}(l)$ is the set of all tree-type diagrams 
 $\tau_k^{(\mu)}$
 constructed with the help of $k$ \mbox{$\mu$-elements} such that 
$|\CE(G(\tau_k^{(\mu)}))| = 2(k-l)+1$ and $|\CV(G(\tau_k^{(\mu)}))|= k-l+2$,
$l=0, \dots, k-1$.

In the third asymptotic regime $(N,c,R)_0^{(\riii)}\to\infty$,
we deduce from (5.5) 
relation
$$
\rcum_k(X^{(\a,3)})= 
\CW_{N,R}(\fT_k^{(\mu,\min)})
(1+o(1)),
\quad (N,c,R)_0^{(\riii)}\to\infty.
$$
In this case 
$|\CE(G(\tau_k^{(\mu)}))|=3$ and  $|\CV(\tau_k^{(\mu)})|=3$ and we get from  (5.10)
equality
 $$
 w^{(\a)}(\tau_k^{(\mu,\min)}) = 
 \int_{-\infty}^{\infty}  
\int_{ -\infty} ^{ \infty}
\ h_{k}^{(\a)}(-x_2)h_{k}^{(\a)}(x_2-x_3)
h_{k}^{(\a)}(x_3)
\, dx_2 \, dx_3.
\eqno (5.14)
$$ 
Having $k$  triangle graphs with oriented labeled edges, we have $6^{k-1}$ 
different minimal tree-type diagrams.
This observation, together with (5.14) proves relation (2.23).
  Theorem 2.3 is proved. $\Box $


\section{Limit theorems for number of walks in random graphs}

Convergence of cumulants of a random variable $\U_n$
makes possible to prove limiting theorems 
for the centered and normalized versions of $\U_n$. 
In particular, 
if there exists
a sequence $(b_n)_{n\ge 1}$ such that
$$
 {1\over b_n} \rcum_k(\U_n) \to \phi_k, \ k\ge 1, \ n\to\infty,
\eqno (6.1)
$$
and $b_n\to\infty$, 
then 
$$
\lim_{n\to\infty, b_n\to\infty } \rcum_k\left( {\U_n- \bE \U_n\over \sqrt{b_n}} \right) = 
\begin{cases}
 \phi_2,  &  \text{if  \   $k=2$} ; \\
0,  \quad   & 
\text {if   \ $k\neq 2$}.
\end{cases}
\eqno (6.2)
$$
The last relation means that the sequence of  centered and normalized random variables 
$
 \chi_n = (\U_n- \bE \U_n)/ \sqrt{b_n} 
$
converges in distribution to a variable $\chi$ with the normal (Gaussian) probability  distribution 
$\CN(0, \phi_2)$,
$$
\chi_n \stackrel{\CL}{\to} \chi, \quad \chi \sim \CN(0, \phi_2), \quad n\to \infty.
\eqno (6.3)
$$

In this section, we prove that convergence (6.2) is valid  for 
normalized random variables $Y_{N,c,R}^{(\a, q)}$ in all of the three asymptotic regimes
determined by Theorems 2.1 and 2.2 while
random variables
$X^{(\a,3)}_{N,c,R}$  converges  in distribution  either to a normal random variable,
when centered and  properly normalized, or to a random variable 
with  Poisson probability distribution, in dependence of the asymptotic behavior
of the parameter $c^3R^2/N^2$ that gives the average number of triangles 
in random graphs considered.   Let us note that the dichotomie between the normal (Gaussian) and the Poisson distribution in the number of sub-graphs
in random graphs is a usual situation, see  \cite{Mau,
P} and references therein.

\vskip 0.5cm 
\subsection{Central Limit Theorems for $Y^{(q)}$-models and $X^{(3)}$-models}
 
 We consider first the case of  $Y^{(q)}$-models.  Let us introduce three random variables 
 $\chi^{(\a,i)}$ that follow centered normal distributions 
 $$
 \chi^{(\a,i)} \sim {\CN}(0, \phi_2^{(\a,i)}), \quad i=1,2,3,
 $$
 with 
 $$
\phi_2^{(\a,1)}=
\begin{cases}
 \Phi_2^{(q,1)}= 2q^2 V_{0}^{2q-1},  &  \text{for  $\a=0$}  , \\
 \Xi_2^{(q,1)} = 2q^2 V_{1}^{2q-2}\,  V_{2},  
 & 
\text {for  $\a=1$,} 
\end{cases}
\eqno (6.4)
$$,
$$
\phi_2^{(\a,2)}=
\begin{cases}
 \Phi_2^{(q,2)}(s),  &  \text{for  $\a=0$}  , \\
 \Xi_2^{(q,2)}(s),    
 & 
\text {for  $\a=1$;} 
\end{cases}
\eqno (6.5)
$$
and 
$$
\phi_2^{(\a,3)}=
\begin{cases}
 \Phi_2^{(q,3)}= 2^{k-1} V_{0},  &  \text{for  $\a=0$}  , \\
 \Xi_2^{(q,1)} = 2^{k-1} V_{kq},   
 & 
\text {for  $\a=1$.} 
\end{cases}
\eqno (6.6)
$$
The right-hand sides of  (6.4), (6.5) and (6.6) are determined with the help of  formulas (2.11), (2.15) and (2.18),
formulas (2.9) and (2.16) and formulas (2.10) and (2.17), respectively.  
 
\vskip 0.2cm 
{\bf Theorem 6.1.} {\it Under conditions of Theorem 2.2, 
 convergence in distribution holds in  three asymptotic regimes   
 %
 of  limiting transition $(N,c,R)_0\to \infty$ (2.6),
 $$
 \chi_{N,c,R}^{(\a,i)}= {1\over 
 \sqrt{b_{N,c,R}^{(i)} }
 }
 \left( Y^{(\a,q)}_{N,c,R} - \bE Y^{(\a,q)}_{N,c,R} \right)
\stackrel{\CL}{\to} 
\chi^{(\a,i)}, \quad i=1,2,3,
\eqno (6.7)
 $$
 respectively, where 
 $$
 b_{N,c,R}^{(i)}= 
 \begin{cases}
 (cR)^{2q-1}/N^{2q-2},  &  \text{if  $cR\gg N$}  , \\
 cR,  
 & 
\text {if  $cR=sN$,} 
\\
cR, & \text{if $cR\ll N$,} 
\end{cases} \quad i=1,2,3.
\eqno (6.8)
 $$
 }
 
\vskip 0.2cm
{\it Proof.} Let us start with the first asymptotic regime of (6.7) when $cR\gg N$. 
Remembering  results of Theorem 2.1 and Theorem 2.2 given by relations (2.8) and (2.15), we 
see that convergence (6.1) is verified with the following choice 
 of variables $b$ and $\U$, 
$$
b = {1\over cR} \quad \hbox{and} \quad  \U= {N^{q-1}\over (cR)^{q-1}} Y^{(\alpha,q)}.
$$
According to (6.2) and (6.3), we get convergence 
$$
{(N/cR)^{q-1} \left(Y^{(\a,q)}- \bE Y^{(\a,q)}\right)\over \sqrt{cR}}
= {Y^{(\a,q)} - \bE Y^{(\a,q)} \over \sqrt{cR(cR/N)^{2q-2}}} 
\stackrel{\CL}{\to} 
 \chi^{(\a,1)}, \quad \hbox{as} \ cR/N\to\infty
$$
with $\chi^{(\a,1)}\sim \CN(0, \phi_2^{(\a,1)})$ and  variance $\phi_2^{(\a,1)}$ given by (6.4).  
This convergence explains the form of the normalizing variable $b_{N,c,R}^{(1)}$ of (6.8). 
The same reasoning shows that convergence (6.2) holds with the choice of 
$
b= cR$ and $\U = Y^{(\a,i)}$ for $i=2$ and $i=3$
and that 
$$
{1\over \sqrt{cR}} \left(Y^{(\a,i)} - \bE Y^{(\a,i)}\right) \stackrel{\CL}{\to} 
\chi^{(\a,i)}, \quad (N,c,R)_0\to\infty, \quad i=2 \ \hbox{and } \ i=3,
$$
in the asymptotic regimes $cR=sN$ and $cR\ll N$, respectively. 
Theorem 6.1 is proved. $\Box$


Similar to Theorem 6.1  statement is valid for $X^{(3)}$-models.
Difference with respect to $Y^{(q)}$-models is that 
 the Central Limit Theorem 
is valid in the third asymptotic regime, but not at the whole extent.
Let us introduce three random variables
$\xi^{(\a,i)}$ that follow centered normal distributions 
 $$
 \xi^{(\a,i)} \sim {\CN}(0, \varphi_2^{(\a,i)}), \quad i=1,2,3,
 $$
with $\varphi_2^{(\a,1)}= \Theta_2^{(\a,\ri)}$  determined by formula (2.25), $\varphi_2^{(\a,2)}= \Theta_2^{(\a,\ri\ri)}(s)$ determined by formula (2.26)
and $\varphi_2^{(\a,3)}= \Theta_2^{(\a, \ri\ri\ri)}$ determined by formulas (2.23) and (2.24). 

\vs 

{\bf Theorem 6.2.} 
{\it Under conditions of Theorem 2.3, 
 convergence in distribution holds in three asymptotic regimes   
 %
 of limiting transition $(N,c,R)_0\to \infty$ (2.6),
 $$
 \xi_{N,c,R}^{(\a,i)}= {1\over 
 \sqrt{d_{N,c,R}^{(i)} }
 }
 \left( X^{(\a,3)} - \bE X^{(\a,3)} \right)
\stackrel{\CL}{\to} 
\xi^{(\a,i)}, \quad i=1,2,3,
 $$
 respectively, where 
 $$
 d_{N,c,R}^{(i)}= 
 \begin{cases}
 {c^5R^3}/N^{4},  &  \text{if  $c^2R\gg N^2$}  , \\
 cR,  
 & 
\text {if  $c^2R=sN^2$,} 
\\
c^{3/2}R/N, & \text{if $c^2R\ll N^2$ and $N^2\ll c^3 R^2$,} 
\end{cases} \quad i=1,2,3.
\eqno (6.9)
 $$
 }

{\it Proof.} In the case of two first asymptotic regimes when $c^2R/N^2 \gg 1$ and $c^2R/N^2=s$, the proof of Theorem 6.2  is based on the use of the results of Theorem 2.3 and relations (6.1) and (6.2). The reasoning repeats the one used to prove Theorem 6.1 and we do not present it here.  

Let us look at the third asymptotic regime when $c^2R/N^2\ll 1$ as $(N,c,R)_0\to\infty$. Convergence (2.23) says that (6.1) is true with the choice of 
$b= c^3R^2/N^2$ and $\U = X^{(\a,3)}$. We see that  if $c^3R^2/N^2\to\infty$, then convergence 
$$
{X^{(\a,3)}- \bE X^{(\a,3)}\over \sqrt{c^3R^2/N^2} } \stackrel{\CL}{\to}  \xi^{(\a,3)}, \quad \xi^{(\a,3)}\sim \CN(0, \varphi_2^{(\a,\ri\ri\ri)}),
\quad (N,c,R)_0\to\infty 
$$
holds. This explains the choice of the normalizing variable $d^{(3)}_{N,c,R}$ of (6.9)  and the last 
condition $c^3R^2\gg N^2$ of (6.9). Theorem 6.2 is proved. $\Box$

\vskip 0.2cm  
Let us briefly discuss results of Theorems 6.1 and 6.2. Observing that 
$$
\bE Y^{(\a,q)}_{N,c,R} = O\big(N (cR/N)^q\big), \quad (N,c,R)_0\to\infty,
$$
we deduce 
 from (6.7) 
that fluctuations of the average number of $q$-step non-closed walks
$\bar Y_{N,c,R} = Y_{N,c,R}/N$ are of the order $O(\bE \bar Y_{N,c,R})/\sqrt {cR}$ in the 
asymptotic regimes when either $cR/N\to\infty$ or $cR/N= O(1)$. 
In contrast to this, in the third asymptotic regime when $cR/N=o(1)$, 
fluctuations of $\bar Y_{N,c,R}$ can be much smaller or much greater
than the average value $\bE \bar Y_{N,c,R}$ in dependence of the ratio between $(cR/N)^q$ and $\sqrt{cR}/N$. In particular,
 in the two-star model determined by $Y^{(\a,2)}$ with $q=2$, the threshold value is given by $cR= N^{2/3}$. 

Regarding $X^{(3)}$-models, we observe that the same is true for 
the random variable $X^{(3)}_{N,c,R}$: its fluctuations are of the 
order $O\big( \bE X^{(3)}_{N,c,R})/\sqrt{cR}$ in the first two asymptotic regimes of Theorem 6.2 while in the third asymptotic regime,
fluctuations of $X^{(3)}_{N,c,R}$ are much smaller than
the mean value $\bE X^{(3)}_{N,c,R}$ only if $c^3R^2/N^2\to \infty$. 
The limiting transition such that  $c^2R/N^2\to\infty$ and $c^3R^2/N^2= O(1)$ will be considered in the next sub-section.

\subsection{Poisson distribution for the number of triangles}

Let us study $X^{(\a, 3)}_{N,c,R}$ (2.4)
in the limiting transition (2.6)
such that
$$
 {c^2R\over N^2}\to 0 \quad {\hbox{and}} \quad 
 {c^3R^2\over N^2}\to \Lambda.
 \eqno  (6.10)
 $$
We denote the limiting transition (6.10) by $(N,c,R)_0^{(\riii)''}\to\infty$.

\vskip 0.2cm

{\bf Theorem 6.3.}
{\it 
If $\a=0$, then 
we have the following convergence in distribution,
$$
T_{N,c,R}^{(0)}= {1\over 6} X^{(0,3)}_{N,c,R}\stackrel{\CL} \to \nu, 
\quad 
(N,c,R)_0^{(\riii)''}\to\infty
\eqno (6.11)
$$
where $\nu$ follows the Poisson probability distribution, 
$$
 \nu \sim \CP(\Lambda  \tilde H^{(0,3)}/6)
\eqno (6.12)
$$
with 
$$
  \tilde H^{(0,3)}= 
{1\over 2\pi} \int_{-\infty}^\infty \big( \tilde h ^{(0)}(p)\big)^3 dp.
$$

\vskip 0.2cm 
If $\a=1$ and $\psi(t)$ of (2.1) is such that 
there exists a random variable $\zeta$ such that $P_\z(x)= P(\z<x)$,
$$
\tilde H^{(1, 3)}_k= 
\int_{-\infty}^\infty s^k dP_\z(s)
\quad {\hbox{and}}
\quad 
\int_{-\infty}^\infty e^{ts} dP_\z(s)<\infty, 
$$
then 
$$
T_{N,c,R}^{(1)}= {1\over 6} X^{(1,3)}_{N,c,R} \stackrel{\CL} \to \nu^{(\z)},
\quad (N,c,R)_0^{(\riii)''}\to\infty, 
\eqno (6.13)
$$
where $\nu^{(\psi)} $ follows the compound Poisson probability distribution 
$
\nu^{(\psi)}\sim \CP(\Lambda/6; P_\z).
$ 
}

{\it Proof.} 
According to (2.23), 
$$
\lim_{(N,c,R)_0^{(\ri\rv)}\to\infty}\rcum_k(tX^{(3)}/6)=
 {\Lambda t^k \over 6} 
\tilde H^{(\a,3)}_k
= {\cal C}_k^{(\a)},
\eqno (6.14)
$$
where $\tilde H_k^{(\a,3)}$ is given by (2.24). 
If $\a=0$, then 
$$
\tilde h_k^{(0)} (p)= \int_{-\infty}^\infty 1^k e^{-\psi^2(x)} e^{-ipx}dx
$$
and $\tilde H_k^{(0,3)}= \CC_k^{(0)}$ do not depend on $k$.
Trivial identity
$$
\sum_{k=1}^\infty  {\CC^{(0)}\over k!} = 
{\Lambda \tilde H^{(0,3)}\over 6 } \, {(e^{t \z }-1)}
$$
shows that $C_k^{(0)}= \CC^{(0)}, k\ge 1$ represent  cumulants of a random variable $\nu$
that follows  
the Poisson probability distribution (6.12). Then (6.11) follows from (6.14) with $\a=0$. 

If  random variable $\z $ exists, then 
the right-hand side of (6.14) 
can be rewritten in the following form
$
\CC_k^{(1)}= {\Lambda t^k} M_k/6, \quad M_k= \int s^k dP_\z(s)= 
\E \z^k.
$
Then 
$$
\sum_{k=1}^\infty  {\CC_k^{(1)}\over k!} = 
{\Lambda \over 6} \left(\E e^{t\z }-1\right)
$$
and (6.13) follows. Theorem 6.3 is proved. $\Box$

%


\subsection{Large number of triangles and finite average  vertex degree }

In this sub-section we  study of $X^{(\a,3)}$ in the 
asymptotic regime  $(N,c,R)_{0}^{(\riii)'}\to\infty$ and 
consider the sequences $c=c(N)$ and $R= R(N)$ such that  
$$
\ {cR\over N}  \to \d
, \quad (N,c,R)^{(\riii)'}_0\to\infty.
\eqno (6.15)
$$
We denote the limiting transition (6.15) by 
$(N,c,R)_{0}^{(\riii)'_a}\to\infty$
and  consider, for simplicity, the case when   
 $R= \d N^\s$ and $c= N^{1-\s}$,  $0<\s <1$.
We are going to  show that, due to Theorem 2.3, 
the total number of triangles in the graph given by 
$T_{N,c,R}= T^{(0)}_{N,c,R}$ (6.11) infinitely increases
in the limit (6.15)
while, according to Theorem 2.1,  
the average vertex degree determined by relation 
$$
 \D_{N,c,R}= {1\over 2N} \sum_{i,j\in \rL_N} a_{ij}^{(N,c,R)} 
= {1\over 2N} Y^{(0,1)}_{N,c,R}
\eqno (6.16)
$$
remains bounded. An explication of this observation 
follows from  elementary equalities,
$$
\bE  \D_{N,c,R}= {cR\over 2N} \cdot {1\over R} \sum_{i,j\in \rL_N} 
\varphi\left({i-j\over R}\right)
=\d  V_0/2(1+o(1)), \quad (N,c,R)^{({\riii})'_a}_0\to\infty,
\eqno (6.17)
$$
and 
$$
\bE T_{N,c,R}= {c^3R^2\over 6N^2}
 \tilde H^{(0,3)} \big(1+o(1)\big)= O(N^{1-\s}), 
\quad (N,c,R)^{({\riii})'_a}_0\to\infty
\eqno (6.18)
$$
that is a consequence of the formula (cf. (2.1))
$$
\bE 
(A^3_{N,c,R})_{ii}=  
\left({c\over N}\right)^3
\sum_{j,l\in \rL_N} 
\varphi\left({i-j\over R}\right) \varphi\left({j-l\over R}\right)
\varphi\left({l-i\over R}\right), 
\eqno (6.19)
$$
where $\varphi(x)= \exp(-\psi^2(x))$. 
Let us formulate the rigorous result.  
\vskip 0.2cm
{\bf Theorem 6.4.} {\it Assume that the infinite family of random variables 
$\{\CA_{N,c,R}, N\in \bN\}$ (2.1)}
{\it with given sequences $c=c(N)$ and $R= R(N)$ are determined on the same probability space. Then under conditions of Theorem 2.1, the following two relations are true, 
$$
P\left(\lim_{ (N,c,R)_0^{(\riii)'_a}\to\infty} 
 \D_{N,c,R}=\d V_0/2\right)=1 \ 
{\hbox{ and }} 
\ P\left(\liminf_{(N,c,R)_0^{(\riii)'_a} \to\infty} T_{N,c,R}= +\infty\right)=1.
\eqno (6.20)
$$
}
\vskip 0.2cm
{\it Proof.} 
Let us start with the last statement of (6.19). Condition (6.15) means that
$c^2R/N^2= \d^2/R\to 0$ and therefore 
 cumulants of $X^{(0,3)}= X^{(0,3)}_{N,c,R}$ verify relation (2.23) 
deduced in the third asymptotic regime of Theorem 2.3,
$$
\lim_{(N,c,R)_0^{(\riii)'_a}\to\infty}\ 
{1\over b_N} \rcum_k(T_{N,c,R})=   1, \quad k=1,2,3,...,
\eqno (6.21)
$$
where  $ b_N = 6\d^2 c\tilde H^{(0,3)}$. 
Let us consider  centered random variables 
$\hat T= T_{N,c,R} - \bE T_{N,c,R}$ and denote their moments
by $\mu_k= \mu_k(\hat T)= \bE (\hat T)^k$. It is known that 
$$
\rcum_1(\hat T)= 0 \quad {\hbox{and}} \quad 
\rcum_k(\hat T)= \rcum_k(T), \ k\ge 2
\eqno (6.22)
$$
and that 
$
\mu_2= \rcum_2(\hat T),
$
$
\mu_3= \rcum_3(\hat T)
$.
The following recurrence is true \cite{BJK-1},
$$
\quad \mu_k= \sum_{l=2}^{k-2}
{{k-1}\choose{l-1}} \, \mu_{k-l} \, \rcum_l(\hat T) + 
\rcum_k(\hat T), \quad k\ge 2.
\eqno (6.23)
$$
Relations (6.20) and (6.21) imply that $\rcum_k(\hat T)= O(b_N)$ and 
therefore
$$
\mu_{2p+1}(\hat T)= O(b_N^p) \quad {\hbox{and}} 
\quad 
\mu_{2p}(\hat T)= O(b_N^p),  \quad (N,c,R)_0^{(\riii)'_a}\to\infty.
\eqno (6.24)
$$
Asymptotic equalities  (6.23) can be proved by recurrence  (6.22). 
Elementary inequality
$$
P\big(| \hat T| \ge  b_N/2\big) \le {\mu_{2p}\over (b_N/2)^{2p}}
$$
combined with (6.24) implies that 
$$
P(\rA_N^{(N,c,R)}) = P\big( \hat T_{N,c,R} < b_N/2\big) \le 4^p 
O\big(b_N^{-p}\big)= O\big(N^{(1-\s)p}\big),
\quad (N,c,R)_0^{(\riii)'_a}\to\infty.
$$
If $p$ is such that $(1-\s)p>1$, then 
$\sum_N P\big(\rA_N^{(N,c,R)}\big)<\infty$ under condition that   (6.15) holds. This convergence  proves the second relation of (6.19). 
\vskip 0.1cm 

Let us consider the centered random variable $\hat \D_N= \D_{N,c,R} - 
\bE \D_{N,c,R}$
and denote 
$$
r_{N,c,R}= \bE \D_{N,c,R} - \d V_0/2.
$$
Then for any $\vep>0$, we have
$$
P\big(| \D_{N,c,R}- \d V_0/2| \ge \vep\big)
\le P\big( | \hat \D_N| \ge \vep-r_{N,c,R}\big)
\le { \rcum_3(\hat \D_N)\over (\vep - r_{N,c,R})^3}
= {cR\over 8 N^3(\vep-r_{N,c,R})^3}.
$$
This implies  the first relation of (6.19). Theorem 6.4 is proved. $\Box$
\vskip 0.1cm 

Let us note that relation (6.17)
shows that the order of the number of triangles in random graphs 
can be arbitrary close to $N$; in particular, the choice of
$R'= \d \ln N$, \mbox{$c'= N/\log N$} still satisfies (6.15). The average vertex degree remains finite and the number of triangles increases as fast as $N/\ln N$ in the ensemble of infinitely increasing random \mbox{graphs (2.1)}.  Theorem 6.4 shows that the ensemble of large Erd\H os-R\'enyi random graphs with large interaction radius solves the random graph collapse problem (see \cite{ACE}). 

To explain results  of Theorem 6.4, let us consider relation (6.19). It says that given $i$, there are in average 
$(cR)^2/N^2$  vertices $j'$ to produce (potentially) a triangle 
with the edge $(i,j')$. Indeed, assuming that the off-spread of  $i$ 
produces  $cR/N$ vertices of the "local world" (see (6.17)), each vertex of this off-spread creates, in its turn, $cR/N$ vertices.  
 This additional edge $(i,j')$ appears with probability approximately $c/N$. Thus, each vertex $i$ participates in $c^3R^2/N^3$ triangles. Summing over $i$ gives the order $c^3R^2/N^2$ (6.18).
 
 One could think also about a large number of $\d$-regular subgraphs (local worlds) with the dilution by $c/N$, $c\to\infty$  that would produce
  an infinitely increasing  number of triangles 
  $N \d^2 \times c/N= \d^2c$, but it seems to be  
 impossible to build $\d$-regular graphs that would have  a dilution of this kind. 
 
 To complete this sub-section, let us consider similar to $T_{N,c,R}$ (6.11) variable
 $$
 S_{N,c,R}= {1\over 8}\  \sum_{ \stackrel{i_1, i_2, i_3, i_4\in L_n}{i_1\neq i_4, \ i_2\neq i_3}} A_{i_1i_2} A_{i_2i_3} 
 A_{i_3i_4} A_{i_4i_1}
 \eqno (6.25)
 $$
that counts the number of squares (quadrangles) in a graph with the adjacency matrix $A_{ij}$. Elementary computation similar to that of (6.18) shows that 
$$
\bE S_{N,c,R} = {c^4R^3\over 8N^3} \tilde H^{(0,4)}(1+o(1)), \quad (N,c,R)_0\to\infty.
\eqno (6.26)
$$
Comparing the right-had side of (6.26) with that of (6.18), we see that 
there exists an asymptotic regime such that the average number of quadrangles  
 in random graphs vanishes while the average number of triangles in this ensemble  remains finite or even increases. 
 With the choice of $c= N^\a, R= N^\b$, $(\a,\b) >0$,  this happens 
 when 
 $$
 3\a + 2\b \ge 3\quad \hbox{and} \quad 4\a + 3\b <1.
 \eqno(6.27)
 $$
 Relations (6.17) shows that the average vertex degree vanishes when (6.27) is true. Let us note that the question of the presence or absence of triangles and/or quadrangles in random graph is closely related 
 with the the sign of a graph  curvature \cite{Bhat,T}. In these studies, the ensemble of Erd\H os-R\'enyi random graphs with large interaction radius can serve as  a useful example due to  possibility to "separate" asymptotic behavior of triangles and quadrangles  
in the limit of infinite  $N,c$ and $R$.


\section{Enumeration of tree-type diagrams}

In Section 4 we have shown that the limiting expressions
for the cumulants of $Y^{(\a,q)}$-model are determined  by
the number of maximal tree-type diagrams
 $\Phi^{(1)}_k= t^{(q)}_k$, in the case of $\a=0$ (4.19). 
 An explicit form of numbers $t^{(q)}_k$ with $q=2$ has been obtained in \cite{K-08},
  with the help of recurrence relations and generating function technique, 
$$
t_k^{(2)}= 2^{2k-1} (k+1)^{k-2}, \quad k\ge 1.
\eqno (7.1)
$$
This sequence, up to the factor $2^{k-1}$, is known in various settings of combinatorial enumeration 
\cite{OEIS} and can be  naturally associated  with the number of trees of $k$ labeled
edges. 

In  \cite{K-2}, it is proved that for general $q\ge 2$, we have 
$$
t_k^{(q)}= 2^{k-1} q^k \big( k(q-1)+1\big)^{k-2}, \quad k\ge 1.
\eqno (7.2)
$$
In this paper, we further generalize  this method of \cite{K-2}
to take into account the multiplicity of edges in the maximal tree-type diagrams 
$\CT_k^{(q)}$ and to obtain an explicit form for $\Xi_k^{(q,1)}$ given by formula (2.18).

\subsection{Pr\"ufer codes for trees and tree-type diagrams  of $k$  elements}

Let us briefly describe the Pr\"ufer codification procedure to get  (7.2) in the case of $q\ge 3$ proposed in \cite{K-2}.
We start with  the case of closed elements $\mu_q$ and construct the color Pr\"ufer code for the
tree-type diagrams $\tau_k^{(q)}\in \fT_k^{(q,\max)}$ (5.3).
We assume that a  tree-type diagram $\tau_k^{(q)}$ 
is constructed with the help of $k$  elements $\mu_q$
labeled by $k$ ordered letters (colors) $\{a,b,\dots,h\}$; 
each edge
of $\mu_q$ is colored in corresponding color. 
The next step is to transform  $\tau_k^{(q)}$ into  a colored diagram 
$\tau_k^{(q,\rc)}$ which will be coded with a Pr\"ufer-type sequence 
$\CP_k^{(q)}$.

To get  $\tau_k^{(q, \rc)}$, we consider $k(q-1)+1$ edges of the
graph 
  $G(\tau_k^{(q)})$ and 
choose among them a root  edge $e_\rho$;
we attribute to it the label  $"0"$;
then we wash out the colors of the edges of $\tau_k^{(q)}$ 
that correspond to 
$e_\rho$. Regarding each of $\mu_q$-elements attached to $e_\rho$, we 
 numerate the edges of  this $\mu_q$-element 
in the 
clockwise direction starting from the 
colorless edge. We say that the ensemble of all $\mu_q$-elements that have  one colorless edge
represent the first layer of $\mu$-elements, we denote this ensemble by 
$\fL_1$. 
Then we consider an element $\mu'_q$ attached to one of the elements of the first layer and say that $\mu'_q$ is the $\mu$-element of the second layer $\fL_2$. 
We erase the color of the edge of $\mu_q'$ attached 
 to the element of the first layer. Then we numerate the color edges of  
 $\mu'_q$ in clockwise direction. 
 Repeat this for all elements of the second layer $\fL_2$.
 Then we pass to $\mu$-elements of the third layer $\fL_3$ and so on. 
  When all color edges are enumerated, we get a new diagram
 $\tau_k^{(q,\rc)}$ ready for the  construction of 
 $\CP_k^{(q)}$.

The color 
 Pr\"ufer sequence 
 $\CP_k^{(q)}$
we are going to construct  
 is given by a sequence of  $k-1$ symbols  taken from the set 
$\fN_k^{(q)}
= \{0, a_1, \dots, a_{q-1}, b_1, \dots, b_{q-1}, \dots, h_{q-1}\}$
of cardinality $|\fN_k^{(q)}|= k(q-1)+1$. 
We consider  $\CP_k^{(q)}$ as  a set of 
 $k-1$ cells (boxes) to fill by recurrence. On the initial step these boxes are empty. 

The recurrent procedure is as follows: in  $\tau_k^{(q,\rc)}$,  
we observe the maximal layer $\fL_m$ of $\mu$-elements; in $\fL_m$, we  find the maximal
element $\mu_q^{\max}$ and consider the element $\tilde \mu_q\in \fL_{m-1}$ that $\mu_q^{\max}$ is attached to; 
we determine the  color $\tilde c$ and the number $ j$ of the edge  
of the element $\tilde \mu_q\in \fL_{m-1}$ that $\mu_q^{\max}$ is attached to; 
put the value $\tilde c_{ j}$ 
into the first cell of the Pr\"ufer-type sequence 
$\CP_k^{(q)}$ and remove the element  $\mu_q^{\max}$ from 
$\tau_k^{(q,\rc)}$. 
As a result of this removal, we get a new diagram $\tau_{k-1}^{(q,\rc)}$.

Now we can repeat 
the procedure described above to fulfill the second cell of 
$\CP_k^{(q)}$.  
When all $k-1$ cells are fulfilled, we get a  color Pr\"ufer
sequence $\CP_k^{(q)}$ that is uniquely determined by $\mu_k^{(q,\rc)}$. 
There is a bijection between the set  of all color 
Pr\"ufer sequences ${\mathfrak P}_k^{(q)}$
and the set of all colored diagrams $\fT_k^{(q,\rc)}$ \cite{K-2}.

By construction, the cardinality of $\fP_k^{(q)}$ is given by
$| \fP_k^{(q)}|= (k(q-1)+1)^{k-1}$. 
Dividing this value by $k(q-1)+1$ that represents the number of possibilities to choose  the root edge $\rho_e$
and multiplying the result by the factor $q^k$ that represents the number of possibilities 
to choose one edge from $\mu$-element whose color is washed out, 
we  get 
the cardinality of the set of all maximal tree-type diagrams, 
$$
|\fT_k^{(q,\max)}|= 2^{k-1}q^{k} (k(q-1)+1)^{k-2},
\eqno (7.5)
$$
where the factor $2^{k-1}$
takes  into  account orientation of  $\mu$-elements of $\tau_k^{(q)}$.
Clearly, the number of tree-type diagrams
$\CT_k^{(q)}$ constructed with the help of 
$k$ $\l$-elements is exactly the same and 
 therefore (7.2) follows from (7.5).

Let us note  that one  can also consider an ordered collection  of   $\l$-elements 
$ \{\l_{r_1}, 
\dots, \l_{r_k}\}$,
with different number of edges
 $1\le r_i\le q$.
 It is not hard to see that the total number of tree-type diagrams constructed with the help of this set  is given by 
$$
t_{r_1,r_2,\dots, r_k} = 2^{k-1} \left(r_1+ \dots + r_k-k+1\right)^{k-2} \, \prod_{i=1}^{k} r_i.
\eqno (7.6)
$$
This relation can be proved with the help of the same color codification procedure as  above.  
With (7.6) in hands, it is easy to understand explicit expression (2.12) and (2.13),
where $\rT^{(q)}(r)$ is the number of all diagrams $d^{(q)}(r)$ obtained 
from one element $\l_q$  by joining its edges in the way such that the number of edges of the graph $E(G(d^{(q)}(r)))=r$
and $\phi_k^{(q,l)}$ is the number of all tree-type diagrams $D_k^{(q)}$
obtained with the help of $k$ $\l_q$-elements
such that the number of edges of the graph $E(G(D_k^{(q)}))=l$.
The numbers $\rT^{(q)}(r)$ are given by known 
 recurrence that we do not present here (see Lemma 6.2 and relation (6.11) of \cite{K-08}).

\subsection{Tree-type diagrams  with multiple edges}

Given a Pr\"ufer sequence $\CP_k^{(q)}$, one can observe that 
there is a symbol $\mathfrak s\in \fN_k^{(q)}$ seen $i$ times in 
$\CP_k^{(q)}$, then the corresponding diagram $\tau_k^{(q,\rc)}$ 
has an edge of multiplicity $i+1$. If there are $s$ different symbols 
seen $i$ times each, then $\tau_k^{(q,\rc)}$ contains $s$ different edges,
each 
of multiplicity $i+1$.

If $ \CP_{k}^{(q)}$ contains
$s_j$  subsets  of $j$ boxes with identical symbols therein, $1\le j\le  k-1$,  
we say that   $ \CP_{k}^{(q)}$  belongs to the equivalence class 
$\bP_k^{(q)}( \s_k) $,  $\s_k= (s_1,\dots, s_{k-1})$.
If $\CP_k^{(q)}\in \fP_k^{(q)}(\s_k)$, then we say that this $\CP_k^{(q)}$ is 
a $\s_k$-Pr\" ufer sequence
and that the corresponding  diagram $\CT_k^{(q)}$ 
is a $\s_k$-tree-type diagram.

 \vskip 0.2cm

{\bf Lemma 7.1.} {\it The number of $ \s_{k}$-type  Pr\" ufer sequences 
  is given by expression
$$
|  \fP_k^{(q)}(\s_{k})|
= {(k-1)!}\, \,  {(k(q-1)+1)!\over (k(q-1)+1- u)!}
  \,  
\
\prod_{i=1}^{k-1} 
{1\over \, (i!)^{s_i} \, s_i!}, 
\eqno (7.7)
$$
where $u= |  \s_{k}|=
s_1+s_2+\cdots + s_{k-1}$.}
 
 \vskip 0.2cm 

{\it Proof of Lemma 7.1.} 
It is known that 
the number 
$$
\CN( \s_{k})={(k-1)!\over   (1!)^{s_1}\,  s_1! \ (2!)^{s_2}\,  s_2!\  \cdots \
((k-1)!)^{s_{k-1}}\,  s_{k-1}!}\ 
\eqno (7.8)
$$
gives the number of possibilities to 
split the set of  $k-1$ elements 
into $u=  |\s_k|$ subsets
of $1$, $2$, \dots, $k-1$ elements, respectively. 
To get a realization of the Pr\"ufer sequence, we have to fill 
$\kappa $ cells by different symbols
taken from the set $\fN_k^{(q)}$ of $k(q-1)+1$ elements. 
This can be done in the following number of ways,
$$
\Upsilon_{k,u}= {(k(q-1)+1)\cdot  k(q-1)\cdots  (k(q-1)-u+1)}
= {(k(q-1)+1)!\over
(k(q-1)-u +1)!}.
\eqno (7.9)
$$
Then the number of  color $\s_k$-Pr\"ufer sequences is given by expression 
$$
|\fP_k^{(q)}(\s_k)| = 
\CN(  \s_{k}) \times \Upsilon_{k,u}= 
{(k-1)!}\ 
  \, u! \,  {{k(q-1)+1}\choose{u}}
\
\prod_{i=1}^{k-1} 
{1\over s_i! \, (i!)^{s_i}}. 
\eqno (7.10)
$$
Lemma 7.1 is proved. 
 $\Box$ 

Using (4.8), we conclude that if $\CT_k^{(q)}$ is a $\s_k$-tree-type diagram, then 
$
w^{(1)}(\CT_k^{(q)})= V_1^{k(q-1)+1- |\s_k|}\,  \prod_{i=1}^{k-1} (V_{i+1})^{s_i}.
$. 
Combining this equality with (7.7) and multiplying the result
by $2^{k-1} q^k(k(q-1)+1)^{-1}$ as  in (7.5), we get 
(2.18).  
\subsection{Proof of relations (2.20), (2.26) and (2.28)}

The second cumulant of the $Y$-model (2.16) is given a total weight of
all tree-type diagrams obtained with the help of two $\l_q$-elements with $q=2$. 
The first term of the right-hand side of (2.20) represents the
total weight of 
maximal tree-type diagrams constructed with the help of two 
$\l_2$-elements.  Using (7.6) with $r_1=2$ and $r_2=2$, we observe that there are 
$t_{2,2}= 8$ such diagrams; each of them has one double edge and two simple edges.
The weight of each diagram is $s^3 V_1^2 V_2$.  
The second term of the right-hand side of (2.20) is represented by 
 diagrams with one simple edge and one triple edge of the weight $s^2 V_1 V_3$; there are $t_{1,2}+t_{2,1}=8$ such diagrams.
 Finally, the third term of the right-hand side of (2.20)  corresponds to diagrams with one quadruple  edge of the weight $sV_4$
and their number is given by $t_{1,1}=2$. Relation (2.20) is proved.

Let us pass to the case of $X^{(3)}$-models. 
The first equality of (2.25) is obvious. 
To prove the second relation of (2.25),
we observe that the limiting expression for the second cumulant of 
$X^{(3)}$-models in the first asymptotic regime is given by 
the total weight of maximal tree-type diagrams constructed with the help of two $\mu_3$-elements. 
According to (2.27), there are  18  such  diagrams; each of them has one double edge and four simple edges. 
According to (5.11), the  weight of each diagram is given by 
$$
w(\tau_2^{(3)})= 
\int_{-\infty}^\infty  
\int_{-\infty}^\infty 
\int_{-\infty}^\infty 
h_1(-x_1) h_1(x_1-x_2) h_2(x_1-x_3) 
h_1(x_2-x_3) h_1(x_3) dx_1 dx_2 dx_3,
$$
where we omitted the superscripts $\a$. 
Then  (2.25) follows.

Relation (2.26) concern the second asymptotic regime of $X^{(3)}$-models. 
The first equality of (2.26) is an obvious consequence of (5.13) considered with $k=1$ and $l=0$. 
The second equality of (2.26) is given by the total weight of 
tree-type diagrams constructed with the help of two $\mu_3$-elements.
The maximal tree-type diagrams produce the weight $\Theta_2^{(\a,\ri)}$ multiplied by $s^2$. There are 6  minimal tree-type diagrams, each of the produces the weight (2.24) multiplied by $s$. This gives the first
term of the second equality of (2.26).

Let us prove relation (2.28). To do this, we   enumerate  tree-type diagrams obtained with the help of
$k$ oriented  elements $\mu_3$. 
Let us consider a partition $\s_{k+1}$ of the set of $k$ triangles
into $u$ subsets, among them there are $s_1$ subsets of one element,
$s_2$ subsets of two elements and so on, $s_k$ subsets of
$k$
elements, such that 
$$
|\s_{k+1}|= s_1+s_2+  \dots + s_k= u \quad \hbox{and} \quad 
\Vert \s_{k+1}\Vert= s_1+ 2 s_2 + \dots + ks_k=k.
\eqno (7.11)
$$
There are 
$$
\CN(\s_{k+1})={k!\over (1!)^{s_1} s_1! \, (2!)^{s_2} s_2! \dots (k!)^{s_k} s_k!}
\eqno  (7.12)
$$
such possibilities under condition (7.11). 
We multiply this expression by the factor
$$
Q_\s= 6^{s_2}\, 6^{2s_3}\dots 6^{(k-1)s_k}= 6^{\Vert \s\Vert - |\s|} 
\eqno (7.13)
$$
that takes into account orientation and root position of $l-1$ elements
in each of $s_l$ subsets.
It remains to construct tree-type diagrams from 
$u$ $\mu$-type elements that gives, according to (7.2), 
$$
t_u^{(3)}= 2^{u-1} 3^u (2u+1)^{u-2}.
\eqno  (7.14)
$$
Combining  expressions (7.2), (7.13) and (7.14) and 
taking into account that each tree 
of $u$ elements produces  the factor $s^u$, 
we get relation (2.28).

\section{On limiting free energy of  $Y$-models}


Returning to definition (1.3), we can write that 
$$
\bE_{{ER}}
\left(e^{- g X_n^{(3)}}\right)
= {1\over (1+e^{-\beta})^{n(n-1)/2} }
\sum_{\g\in \G_n}
 \exp\left\{ - {\beta\over 2} \Tr A^2(\g) - g \Tr A^{3}(\g) \right\}, \ \beta = \ln (n/c),
\eqno (8.1)
$$
where $\G_n$ denotes the family of all possible simple non-oriented loop-less graphs on $n$ vertices
and $A(\g)$ is the adjacency matrix of the graph $\g \in \G_n$ (see \cite{K-08} for more details). 
This means that variable $Z_n^{(X^{(3)})}(g)$ 
(1.3) can be determined, after extraction of evident simple term 
$n(n-1)/2 \ln(1- e^{-\beta})$,
 as the logarithm of the normalization constant $\CC_n(\beta,g)$ of the probability distribution  such that 
$$
P(\g) = {1\over \CC_n(\beta,g)}
 \exp\left\{ - {\beta\over 2} \Tr A^2(\g) - g \Tr A^{3}(\g) \right\}.
\eqno (8.2)
$$
Therefore, results (1.5) on the asymptotic properties of terms of the cumulant expansion of $Z_n^{(X^{(3)})}(g)$
can be used in the study of 
 analyticity domain of this variable and thus are important with respect to the 
 presence or absence of  phase transitions in  discrete random matrix model determined by (8.1). 
 This  problem has been  discussed  in a number of papers on one or another level of rigor (see \cite{M,NRS,RY,Yin}). 
  Probability distribution (8.2) known as the one of the exponential random graphs  attracted much interest, see 
 \cite{HL,Jon} for initial studies. Probability distribution modified 
by the presence of 
$q$-step walk in random graphs, and thus related with $Y$-models have been considered in papers \cite{AC,PN-1}.

In this connection, results of Theorems 2.1 and 2.2 
on the asymptotic behavior of 
 terms of formal cumulant expansion
$$
Z_{N,c,R}(g) = \log \bE_{N,c,R} \big(g Y^{(\a,q)} \big)= 
\sum_{k\ge 1}
{g^k \rcum_k\big( Y^{(\a,q)}\big)/ k!},
\eqno (8.3)
$$
can be also interpreted as the ones related with  asymptotic behavior of the 
free energy of a discrete matrix model. Indeed, we can write that 
$$
\log  \bE_{N,c,R} \big(g Y^{(\a,q)} \big)= 
\log \CC_{N,c,R}^{(\a,q)}(g) - \sum_{-n\le i< j\le n}
\log \left( 1 + {c\over N} e^{-\psi^2((i-j)/R)}\right)
$$
where 
$$\CC_{N,c,R}^{(\a,q)}(g) = \sum_{\g \in \G_N}  
\exp\left\{ - \sum_{i<j} 
\left(\b + \psi^2\left({i-j\over R}\right)\right) 
A_{ij}(\g)
+ g Y^{(\a,q)}(\g) \right\}.
\eqno (8.4)
$$
Similarly to  (8.1), one can say that 
  (8.4) represents  
normalization constant 
of the probability distribution of the exponential random graphs model. 
 The free energy of the matrix  model (8.4) 
 is given by 
 $$
 \CF_{N,c,R}(g) = 
{1\over N} \log \CC_{N,c,R}(g)
= {1\over N} \bE_{N,c,R} \big(g Y^{(\a,q)} \big) + 
{cR\over N} V_0\big( 1+ o(1)\big), \quad (N,c,R)_0\to\infty.
\eqno (8.5)
$$
Then we can rewrite results of Theorems 2.1 and 2.2 as follows,
$$
\CF_{N,c,R}^{(\a,q)}(g)= 
{cR\over N}\CF^{(\a,q)}_i(g)\big( 1+ o(1)\big),
\quad {(N,c,R)_0^{(i)}\to\infty},
\quad i=1,2,3,
\eqno (8.6)
$$
where 
$$
\quad   \CF^{(\a,q)}_i(g)
= \sum_{k\ge 1} {g^k\over k!} \CF_k^{(\a,q)}(g),
\eqno (8.7)
$$
and where, according to relations (2.8), (2.9) and (2.10),
$$
\CF_k^{(0,q)}(g) = 
\begin{cases}
 \Phi_k^{(q,1)},  &  \text{in the limit  $(N,c,R)_0^{(1)}\to\infty$} ; \\
\Phi_k^{(q,2)}(s),     & 
\text {in the limit $(N,c,R)_0^{(2)}\to\infty$};
\\
\Phi_k^{(q,3)}, & \text{in the limit $(N,c,R)_0^{(3)}\to\infty$}.
\end{cases}
\eqno (8.8)
$$
When passing from (8.5) to (8.6), we have omitted the term 
$cRV_0/N(1+o(1))$ that does not depend on $g$. Similar to (8.8) relations can be written in the case of $\a=1$. 
Finally, let us note that in the first asymptotic regime, 
we  consider (8.3) with normalized variable 
$\hat  Y_{N,c,R}= (N /cR)^{q-1} Y_{N,c,R}$  instead of $Y_{N,c,R}$
(see (4.18)).


\subsection{$Y^{(0,q)}$-models in three asymptotic regimes }

Remembering  that the average value of vertex degree 
of the random graphs (2.1)
is of the order $\d=cR/N$ (6.16),  
we can say that the first limiting transition 
of Theorem 2.1 and Theorem 2.2 given by  $\delta\gg 1$ corresponds to the asymptotic regime
of dense  graphs. 
In this case, rewriting relations (2.11) for $\Phi_k^{(q,1)}$
in the form
$$
\Phi_k^{(q,1)}= \left( 2q(q-1)V_0^{q-1}\right)^k \cdot {V_0\over 2 \big(k(q-1)+1\big)^2} \left( k+ {1\over q-1}\right)^k,
\eqno (8.9)
$$
we conclude that the value
$
g_0= 1/(  2eq(q-1) V_0^{q-1})
$
is the critical one for the convergence of the series (8.1). 
In the second asymptotic regime of sparse random graphs with  
$\d =  O(1)$, one can use the estimate from below of the form (8.9),
$
\Phi_k^{(q,2)}(s) \ge 
s^{k(q-1)+1} \Phi_k^{(q,2)}
$
and   conclude that if $g>g_0$, then the infinite series of (8.7) diverges for any given $s$. 
Then one can say that 
 in the case of sparse random graphs, the phase transition exists in the domain $\{(s,g): \ s\in \bR, g\le g_0\}$. 
 In the third asymptotic regime of very sparse random sgraphs, $\d \ll 1$, 
 the limiting expression 
$\CF_3^{(\a,q)}(g)$ exists and is analytical for all $g$. In this asymptotic regime, the exponential graph model (8.4) should not 
exhibit any phase transition behavior. 
These three observations affirm and generalize the statements of the presence or absence of phase transitions obtained in the case of 
two-star model, $Y^{(q)}$ with $q=2$  \cite{AC,BJK,CD,M,RY}.

\subsection{$Y^{(1,q)}$-models in dense graph regime}

 
 In this subsection we use the Bernoulli and Poisson random variables 
 to get the upper estimates for the cumulants of $Y$-models. 
   Let us consider the sum of $k(q-1)+1$ i.i.d. Bernoulli random variables 
  $\xi_i$ 
  and write down its $(k-1)$-th moment in the following form,
$$
\E \left( \xi_1 + \xi_2+ \dots + \xi_{k(q-1)+1}\right)^{k-1}=
\sum_{1\le i_1, i_2, \dots, i_{k-1} \le k(q-1)+1 } \E
 \big( \xi_{i_1}  \xi_{i_2} \cdots 
\xi_{i_{k}}\big).
$$
Considering 
a realization  $\la i_1,\dots, i_{k-1}\ra $ such that there are 
$s_j$ subsets of $j$ identical  elements,
we can write that 
$\E  \big( \xi_{i_1}  \xi_{i_2} \cdots 
\xi_{i_{k-1}}\big)=
p^{u}$, where  $u=s_1+ s_2+ \dots + s_{k-1}$.  
Then 
 $$
\E \left( \xi_1 + \xi_2+ \dots + \xi_{k(q-1)+1}\right)^{k-1}=
\sum_{u=1}^{k-1} p^u \,  \Upsilon_{k,u}
\sum_{ \stackrel{ \s_{k}= (s_1, \dots, s_{k-1})} 
{ | \s_{k}| = u, \, 
\Vert  \s_{k}\Vert = k-1}} \ 
\CN(  \s_{k}) 
$$
$$=
\sum_{u=1}^{k-1} 
{(k-1)!}\ 
  \, u! \,  {{k(q-1)+1}\choose{u}}
\sum_{ \stackrel{ \s_k= (s_1, \dots, s_{k-1})} 
{ |\s_{k}| = u, \, 
\Vert  \s_{k}\Vert = k-1}} \ \prod_{j=1}^{k-1} 
{p^{s_j}\over s_j! \, (i!)^{s_j}}, 
\eqno (8.11)
$$
where $\Upsilon_{k,u}$ and $\CN(  \s_{k})$ are determined by (7.9) and 
(7.8), respectively and where  
 the sum runs over $\s_{k}$ of (2.19).  
Regarding (2.18) with $V_{i+1}= 1$ and using (8.11),
we can write that 
$$
\Phi_k^{(q,1)} 
= 
{q^k\, (k-1)!\over k(q-1)+1} \, 
 \ \sum_{u=1}^{k-1}  \, {(k(q-1)+1)\cdots (k(q-1)-u+2)\over p^u}
 \ 
\sum_{ \stackrel{ \s_k= (s_1, \dots, s_{k-1})} 
{ | \s_k| = u, \, 
\Vert  \s_k\Vert = k-1}} \ 
\prod_{i=1}^{k-1} {1\over s_i!} \left( {p\over  i! }\right)^{s_i}
$$
$$
\le 
{q^k\over p^{k-1} \big(k(q-1)+1\big)} \, 
 \E 
 \left(
 \xi_1 + \dots +  \xi_{k(q-1)+1}
 \right)^{k-1}.
 $$
Taking $p=1/(q-1)$, we get inequality 
$
\Phi_k^{(q,1)} 
\le  q^{2k-1}\, 
 \E \big(
\hat  \xi_1 + \dots +  \hat \xi_{k(q-1)+1}\big)^{k-1}/(k(q-1)+1),
 $
 where 
$$
\hat \xi_i = 
\begin{cases}
1,  &  \text{with probability \   $1/(q-1)$} , \\
0,  \quad   & 
\text {with probability    \ $1 - 1/(q-1)$}.
\end{cases}
\eqno (8.12)
$$
 Let us  consider cumulants of $Y^{(1,q)}$-model and rewrite the right-hand side of (2.18) as  follows,
$$
\Xi^{(q,1)}_k = 
{q^k\, (k-1)!\over k(q-1)+1} {V_1}^{k(q-1)+1}   
 \ \sum_{m=1}^{k-1}  \, m! \,  {{k(q-1)+1}\choose{m}}\,  
 \sum_{ \stackrel{(s_1, \dots, s_{k-1})} 
{ |\s_k| = m, \, 
\Vert \s_k\Vert = k-1}}   
\prod_{i=1}^{k-1} {1\over s_i!} 
\left( {  V_{i+1}/V_1\over i!}\right)^{s_i}
$$
Remembering that  
  $V_{i+1}/V_1$ represents the $i$-th moment of a random variable $\kappa$,
$$
V_{i+1}/V_1= \int (1+\psi^2(x))^{i+1} e^{-\psi^2(x)}dx
/\int (1+\psi^2(x)) e^{-\psi^2(x)} dx= \int (1+\psi^2(x))^i f(x)dx,
$$
where $\kappa = 1+ \psi^2(\zeta)$ and $\zeta $ has a probability  density
$
f(x) = (1+\psi^2(x) ) e^{-\psi^2(x)}/V_1,
$
we get the following representation of the the limiting cumulants 
$\Xi_k^{(q,1)}$,
$$
\Xi_k^{(q,1)}\le 
{q^{2k-1}\, {V_1}^{k(q-1)+1}\over k(q-1)+1} 
\E\left( \kappa_1 \hat \xi_1
+\dots +
\kappa_{k(q-1)+1} \hat \xi_{k(q-1)+1}
\right)^{k-1},
\eqno (8.13)
$$
with random variables $\hat \xi_i$ determined by (8.12).
We see that 
the normalized limiting cumulants 
$\Xi_k^{(q,1)}/q^{2} V_1^{q-1}$ admit asymptotic upper bound by $(k-1)$-th moment
of the random variable
$$
\Lambda_k^{(q)}= \kappa_1 \hat \xi_1
+\dots +
\kappa_{k(q-1)+1} \hat \xi_{k(q-1)+1},
$$
 If $q\to\infty$, then 
$\Lambda_k^{(q)}$ converges in law to a random variable 
$\Lambda_{k\ep}$ that follows the compound Poisson distribution $\CP^{(\kappa)}(k\ep)$ with mean value  $k\ep $, $\ep = V_2/V_1$.  
$$
\Lambda_k^{(q)}\stackrel{\CL}{\to}  \Lambda_{k\ep },\quad q\to \infty, \quad
\Lambda_k\sim \CP^{(\ep )}(k\ep  ).
\eqno (8.14)
$$
Asymptotic properties  of  moments of the compound Poisson distribution,
\mbox{$
\CM_{k-1} = \E \Lambda_{k\ep }^{k-1}$}  in the limit  $k\to\infty$
have been studied in papers \cite{K-1,K-2}. 
We have 
$$
\CM_{k-1} = \left( k\ep  e^{\upsilon(s )}(1+o(1)\right)^k,
\quad k\to\infty,
\eqno (8.15)
$$
where 
$$\CS(x)= \sum_{k\ge 0}  { x^l\over l!} \, \E \kappa^l\, ,\quad  
\upsilon(s)= {\CS(u  ) -1\over  \CS'(u)} - 1 + \ln \CS'(u)
$$
and 
the value of $u$ is such that 
$
u\CS'(u)= {1/ s}.
$

Using upper estimate  (8.13), convergence (8.14)  and 
asymptotic relation  (8.15), we can put forward a conjecture that  
the following upper bound for the limiting cumulants $\Xi^{(q,1)}_k$
$$
 \Xi^{(q,1)} \le 
{q^{2k-1} V_1^{k(q-1)+1}\over k(q-1)+1} \left( k \ep e^{\upsilon(s)}
(1+o(1)\right)^k
\eqno (8.16)
$$
holds for large values of $q$ and $k$. 
Using elementary computations based on (8.16), we can 
argue that the exponential graph model (8.4)
with $Y^{(1,q)}$ replaced by 
$$
\tilde  Y^{(1,q)}(\g)= {1\over q^2 V_1^{q-1}} \left({N\over cR}\right)^{q-1} Y^{(1, q)}(\gamma)
$$
considered for large values of $q$ in the asymptotic regime of dense graphs, 
might have the free energy per cite analytical for all 
$
t< t_1 = ( V_2 \min_{s>0}e^{\upsilon(s) +1})^{-1}.
$
This means that the model (8.4) should not have  phase transitions
for $t<t_1$.

\end{document}